\documentclass[]{amsart}
\usepackage[lmargin=3cm,tmargin=3.5cm,bmargin=3.5cm,rmargin=3cm]{geometry} % margens
\usepackage[table]{xcolor}
\usepackage{tabularx}
\usepackage{mathrsfs,amsmath, amsfonts, amssymb, amscd}
\usepackage{tikz}
\usepackage{graphicx,graphics,epsfig}
\usepackage{epstopdf}
\usepackage{hyperref}
\usepackage{marginnote}
\usepackage{enumerate}

\newtheorem{theorem}{Theorem}[section]
\newtheorem{proposition}[theorem]{Proposition}
\newtheorem{corollary}[theorem]{Corollary}
\newtheorem{lemma}[theorem]{Lemma}
\newtheorem{remark}[theorem]{Remark}

\newcommand{\p}{\partial}

\newcommand{\R}{{\mathbb R}}

\newcommand{\C}{{\mathbb C}}

\title[$H^1$-Critical Schr\"odinger-Debye System]{On the Unboundedness of Higher Regularity Sobolev Norms of Solutions for the Critical Schr\"odinger-Debye System with Vanishing Relaxation Delay}

\author{Ad\'an  J. Corcho}

\author{Jorge Drumond Silva}

\address{\textbf{Ad\'an J. Corcho}
\newline
Instituto de Matem\'atica.
\newline
Universidade Federal do Rio de Janeiro-UFRJ.
\newline
Ilha do Fund\~{a}o, 21945-970. Rio de Janeiro-RJ, Brazil.
\newline
Rio de Janeiro-RJ, Brazil.} \email{adan@im.ufrj.br}

\address{\textbf{Jorge Drumond Silva}
\newline
Center for Mathematical Analysis, Geometry and Dynamical Systems,
\newline
Department of Mathematics,
\newline
Instituto Superior T\'ecnico, Universidade de Lisboa
\newline
Av. Rovisco Pais, 1049-001 Lisboa, Portugal.}
\email{jsilva@math.tecnico.ulisboa.pt }

\thanks{A. J. Corcho was partially supported by CAPES and CNPq (Edital Universal - 481715/2012-6), Brazil}
\thanks{J. Drumond Silva was partially supported by FCT/Portugal through UID/MAT/04459/2013}


\date{\today}

\begin{document}

\setcounter{page}{1}

\begin{abstract}
	We consider the Schr\"odinger-Debye system in $\R^n$, for $n=3,4$. Developing on previously known local well-posedness results, we start by establishing global well-posedness in $H^1({\R}^3)\times L^2({\R}^3)$ for a broad  class of initial data. We then concentrate on the initial value problem in $n=4$, which is the energy-critical dimension for the corresponding cubic nonlinear Schr\"odinger equation. We start by proving local well-posedness in  $H^1({\R}^4)\times H^1({\R}^4)$. Then, for the focusing case of the system, we derive a virial type identity and use it to prove that for radially symmetric smooth initial data with negative energy, there is a positive time $T_0$, depending only on the data, for which,  either the $H^1({\R}^4)\times H^1({\R}^4)$ solutions blow-up in $[0,T_0]$, or the higher regularity Sobolev norms are unbounded on the intervals $[0, T]$, for $T>T_0$, as the delay parameter vanishes. We finish by presenting a global well-posedness result for
	regular initial data which is small in the $H^1({\R}^4)\times H^1({\R}^4)$ norm.
\end{abstract}

\maketitle

\section{\bf{Introduction}.} 
The purpose of this paper is to present new results concerning the dynamics of the Cauchy problem associated to the Schr\"odinger-Debye system, for spatial dimensions three ($n=3$) and four ($n=4$). More precisely,  we consider the system given by the coupled equations:
\begin{equation}\label{SD}
\begin{cases}
iu_t+\frac{1}{2}\Delta u=uv, & (x,t)\in \mathbb{R}^{n}\times \mathbb{R},\\
\mu v_t + v = \lambda|u|^{2}, & \mu >0,\, \lambda=\pm 1, \\
u(x,0)=u_{0}(x),\quad v(x,0)=v_{0}(x), &
\end{cases}
\end{equation}
where  $\Delta=\sum\limits_{j=1}^{n}\p^2_{x_j}$ is the Laplacian operator on $\R^n$, $u=u(x,t)$ is a complex-valued function and $v=v(x,t)$ is a real-valued
function. This system models the propagation of an electromagnetic wave through a nonresonant medium, whose
nonlinear polarization lags behind the induced electric field (see \cite{Newell} for more physical details).
We notice that in the absence of delay  ($\mu=0$), representing an instantaneous polarization response, the 
system (\ref{SD}) reduces to the cubic non-linear
 Schr\"odinger equation (cubic NLS):
\begin{equation}\label{CNLS}
\begin{cases}
iu_t+ \tfrac{1}{2}\Delta u = \lambda |u|^2u,\quad (x,t)\in \mathbb{R}^n\times \mathbb{R},\\
u(x,0)=u_0(x).
\end{cases}
\end{equation}
The cases $\lambda =-1$ and $\lambda =1$  model focusing and defocusing nonlinearities, respectively. We classify the coupling in \eqref{SD} analogously.

\smallskip
In 1999 Fibich and Papanicolau (\cite{Fibich-Papanicolau}) used an extension of an adiabatic approach (developed earlier by Fibich for \eqref{CNLS}) to a general modulation theory in order to study the formation of singularities of self-focusing solutions for small perturbations of the cubic NLS equation (\ref{CNLS}), in the $L^2$- scaling critical dimension; that is, when $n=2$ and $\lambda =-1$. Among the examples of such perturbations to which this technique was applied, they considered, for instance, dispersive
saturating nonlinearities,
self-focusing with Debye relaxa\-tion, the Davey-Stewartson equations, self-focusing in optical fiber arrays 
and the effect of randomness. 
In the particular case of the perturbation of the cubic NLS modeled by the Schr\"odinger-Debye system (\ref{SD}) for $n=2$, the authors of \cite{Fibich-Papanicolau} addressed the question of whether Debye relaxation can arrest self-focusing when  $0< \mu \ll 1$. As a result of this study, it was concluded that self-focusing becomes temporally asymmetrical and thus the modulation theory cannot be conclusive regarding the formation of singularities. On the other hand, from a numerical approach, Besse and Bid\'egaray (\cite{Besse-Bidegaray}) used two different methods suggesting the blow-up, in finite time, of the  $L^{\infty}$-norm for solution $u$ for the specific initial data $u_0(x,y)=\text{e}^{-(x^2+y^2)}$ and $v_0=-|u_0|^2$. 
Recently, the above question was answered in \cite{Corcho-Oliveira-Silva}, where it was proved that,
in the two dimensional case ($n=2$), singularities do not form in finite time, for initial data $(u_0, v_0)$ belonging to the space $H^1(\R^2)\times L^2(\R^2)$.

\smallskip
In this paper, we will show that, for the focusing ($\lambda=-1$) case of system (\ref{SD}) and small relaxation parameter $\mu$, in dimension $n=4$, i.e., for the related $H^1$-critical dimension of the cubic NLS model (\ref{CNLS}), the solutions corresponding to negative energy radially symmetric smooth initial data either blow-up in finite time or their higher regularity Sobolev norms over any large enough time interval are unbounded, as 
$\mu \to 0$.

\smallskip
Before establishing the main results we will review some important properties of the solutions 
to the system (\ref{SD}). The flow preserves the $L^2$-norm of the solution $u$, that is,
\begin{equation}\label{Conservation Law}
\int_{\R^n}|u(x,t)|^{2}dx=\int_{\R^n}|u_{0}
(x)|^{2}dx.
\end{equation}
Also, the following pseudo-Hamiltonian structure holds:
\begin{equation}\label{Energy-1}
\frac{d}{dt}E(t)=2\lambda\mu\int_{\R^n}(v_t(x,t))^{2}dx,
\end{equation}
where
\begin{equation}\label{Energy-2}
E(t)=\int_{\R^n}\Bigl(|\nabla u|^{2}+\lambda |u|^{4}-\lambda \mu^{2}(v_t)^2\Bigl)dx
=\int_{\R^n}\Bigl(|\nabla u|^{2}+2v|u|^2 -\lambda v^2\Bigl)dx.
\end{equation}
This energy integral is well defined as long as $v \in L^2(\R^n)$ and
 $u \in H^1(\R^n)$, with the
Sobolev embedding theorem permitting the $L^4$ norm of $u$ to be controlled by $H^1$, i.e. for $n\leq 4$
 (corresponding to the $H^1$-subcritical and critical dimensions of the cubic NLS). 
Note  from \eqref{Energy-1} that this pseudo-Hamiltonian is not conserved. Although we can immediately infer its monotonicity, depending on the sign of $\lambda$: increases in time, when $\lambda=1$, or decreases, when $\lambda=-1$.

\smallskip
The system (\ref{SD}) can be decoupled  by solving the second equation with respect to $v$,
\begin{equation}\label{equation-v}
v(t)= e^{-t/\mu}v_0(x)+ \tfrac{\lambda}{\mu}\int_0^t\,
e^{-(t-\tau)/\mu}|u(\tau)|^2\,d\tau,
\end{equation}
to obtain the integro-differential equation for $u$,
\begin{equation}\label{SD-IDF}
\begin{cases}
iu_t+\tfrac{1}{2}\Delta u = e^{-t/\mu}uv_0(x)
+\tfrac{\lambda}{\mu}u\displaystyle \int_0^t\,e^{-(t-\tau)/\mu}|u(\tau)|^2d\tau,&  x\in {\R}^n,\; t\ge 0, \\
u(x,0) = u_0(x).
\end{cases}
\end{equation}
Heuristically, \eqref{SD-IDF} illustrates the 
property that, as time increases, the system steadily tends towards resembling a cubic NLS equation, with the
speed of that approximation increasing as $\mu$ decreases, due to its 
effect on the negative exponentials. Physically,
of course, this  reflects the Debye polarization delay, which decreases with $\mu$. 

\smallskip
Another instance of this phenomenon is obtained by applying $L^1$ norms to \eqref{equation-v}, by which we can conclude that the bound for $\|v(\cdot,t)\|_{L^1(\R^n)}$ shifts from
$\|v_0\|_{L^1(\R^n)}$, at $t=0$, to $\|u_0\|_{L^2(\R^n)}^2$, as $t \to +\infty$,
\begin{eqnarray*}
	\|v(\cdot,t)\|_{L^1(\R^n)} &\leq& e^{-t/\mu}\|v_0\|_{L^1(\R^n)} + 
	\tfrac{1}{\mu}\int_0^t\,
	e^{-(t-\tau)/\mu}\int_{\R^n}|u(x,\tau)|^2\,dx\,d\tau,\\
	&=&e^{-t/\mu}\|v_0\|_{L^1(\R^n)} + (1-e^{-t/\mu})\|u_0\|_{L^2(\R^n)}^2,
\end{eqnarray*}
yielding the following a priori bound for the $L^1$ norm of $v$, independently of the dimension or
the sign of $\lambda$, for the whole time interval of existence of the solution,
\begin{equation} \label{apriori_v}
	\|v(\cdot,t)\|_{L^1(\R^n)}\leq \|u_0\|_{L^2(\R^n)}^2 + e^{-t/\mu}\big(\|v_0\|_{L^1(\R^n)}- \|u_0\|_{L^2(\R^n)}^2\big)
	 \leq \|v_0\|_{L^1(\R^n)} + \|u_0\|_{L^2(\R^n)}^2.
\end{equation} 

\smallskip
Unlike  the cubic non-linear Schr\"odinger equation (\ref{CNLS}), the solutions of (\ref{SD}) are not invariant under scaling, but the Debye relaxation allows a dynamical rescaling over the delay parameter $\mu$. Indeed,  if $(u,v)$ is a solution to (\ref{SD}) for a value of $\mu >0$, then
\begin{equation}\label{parameter-rescaling}
\Bigl(\tilde{u}(x,t),\tilde{v}(x,t)\Bigl)=\left(\mu^{1/2}u(\mu^{1/2}x, \mu t),\,\mu v(\mu^{1/2}x, \mu t)\right)
\end{equation}
yields a solution to (\ref{SD}) for $\mu =1$, with initial data $\bigl(\tilde{u}_0(x),\tilde{v}_0(x)\bigl)=\left(\mu^{1/2}u_0(\mu^{1/2}x),\,\mu v_0(\mu^{1/2}x)\right)$. Then, as was already remarked in \cite{Besse-Bidegaray}, we see that the formation of singularities, in case they occur, does not depend on the size of $\mu$, as long as this parameter stays positive, although for necessarily different initial data, according to the previous scaling formula.

\subsection{Overview of known well-posedness results}
Many results, concerning local well-posedness for the Cauchy problem (\ref{SD}) with initial data $(u_{0},v_{0})$ in Sobolev spaces $H^{s}(\R^n)\times H^{\kappa}(\R^n)$, $1\le n\le 3$, have been obtained by applying a fixed-point procedure to the Duhamel formulation asso\-ciated to the integro-differential equation (\ref{SD-IDF}), combined with classical smoothing effects for the  Schr\"odinger group $e^{it\Delta/2}$. We refer to the works \cite{Bidegaray1, Bidegaray2, Corcho-Linares} for more details. Recently (see \cite{Corcho-Matheus, Corcho-Oliveira-Silva}), more general results about local and global well-posedness were obtained in the framework of Bourgain's spaces, by using a fixed-point procedure applied directly to the integral Duhamel formulation for the system (\ref{SD}) itself. These latest results 
contain the previous ones in \cite{Bidegaray1, Bidegaray2, Corcho-Linares} as particular cases. We summarize them as follows:

\smallskip
\begin{theorem}[\cite{Corcho-Matheus, Corcho-Oliveira-Silva}]\label{LWP-BourgainSpaces}
Let $n=1,2,3$. Then, for any $(u_0,v_0)\in H^s(\R^n) \times H^{\kappa}(\R^n)$, with $s$ and $\kappa$ satisfying the conditions:
\begin{enumerate}
\item [(a)] $|s|-\frac12\le \kappa < \min\bigl\{s+\frac12,\; 2s+\frac12\bigl\}$\; and\; $s> -\frac14$\; for $n=1$\,
\emph{(}see \cite{Corcho-Matheus}\emph{)},

\smallskip
\item [(b)]$\max\{0, s-1\} \le \kappa \le \min\{2s, s+1\}$ for $n=2, 3$\,
\emph{(}see \cite{Corcho-Oliveira-Silva}\emph{)},
\end{enumerate}
there exists a time $T=T(\|u_0\|_{H^{s}}, \|v_0\|_{H^{\kappa}})>0$ and
a unique solution $(u(t),v(t))$ of the initial value problem
(\ref{SD}) in the time interval $[0,T]$, satisfying
$$(u, v)\in C\left([0,T]; H^s(\R) \times H^{\kappa}(\R)\right).$$
Moreover, the map $(u_0,v_0) \longmapsto (u(t),v(t))$ is locally
Lipschitz. In addition, when $-3/14< s=\kappa \le 0$, for $n=1$, and $(s, \kappa)=(1,0)$, for $n=2$, the
local solutions can be extended to any time interval $[0,T]$.
\end{theorem}

\smallskip
Figures \ref{Region-A} and \ref{Region-B} represent the regions $\mathcal{W}_1$ and $\mathcal{W}_{2,3}$
in the $(s, \kappa)$ plane, 
corresponding to the sets of Sobolev indices for which local well-posedness (l.w.p.) has been  
established for $n=1$, in \cite{Corcho-Matheus}, and for $n=2,3$, in \cite{Corcho-Oliveira-Silva}, 
as described in Theorem \ref{LWP-BourgainSpaces}.
\begin{figure}[htp]
\begin{minipage}[b]{0.45\linewidth}
\centering
\begin{tikzpicture}
%\draw[very thin,color=gray] (-1.9,-1.9) grid (3.9,3.9);
\draw[->] (-2,0)--(4,0) node[below] {$s$};
\draw[->] (0,-1.5)--(0,4) node[right] {$\kappa$};
\filldraw[color=yellow!50]
(-0.25,-0.25)--(-0.25,0)--(0,0.5)--(3.5,4)--(4,4)--(4,3.5)--(0,-0.5)--(-0.25,-0.25);
\draw[thick, dashed](-0.25,-0.25)--(-0.25,0)--(0,0.5)--(3.5,4);
\draw[thick](4,3.5)--(0,-0.5)--(-0.25,-0.25);
\draw[very thin](-0.25,0)--(0.5,0);
\draw[very thin] (0,-0.5)--(0,0.5);
\node at (2,2){\small{$\boldsymbol{\mathcal{W}_1}$}};
\node at (-0.9,0.4){\small{$\kappa=2s+\frac12$}};
\node at (1.2,2.7){\small{$\kappa=s+\frac12$}};
\node at (1.1,-0.4){\small{$\kappa=|s|-\frac12$}};
\end{tikzpicture}
\caption{l.w.p. for $n=1$ (\cite{Corcho-Matheus})}\label{Region-A}
\end{minipage}\hfill
\begin{minipage}[b]{0.45\linewidth}
\vspace{0.7cm} \centering
\begin{tikzpicture}
\draw[->] (-2,0)--(4,0) node[below] {$s$};
\draw[->] (0,-1.5)--(0,4) node[right] {$\kappa$};
\filldraw[color=yellow!50](0,0)--(1,0)--(4,3)--(4,4)--(3,4)--(1,2)--(0,0);
\draw[thick](3,4)--(1,2)--(0,0)--(1,0)--(4,3);
\node at (0,0){$\bullet$};
\node at (1,0){$\bullet$};
\node at (2,2){$\boldsymbol{\mathcal{W}_{2,3}}$};
\node at (1.5,3.5){\small{$\kappa=s+1$}};
\node at (3.2,1.2){\small{$\kappa=s-1$}};
\node at (0.2,1.5){\small{$\kappa=2s$}};
\end{tikzpicture}
\caption{l.w.p. for $n=2,3$ (\cite{Corcho-Oliveira-Silva})}\label{Region-B}
\end{minipage}
\end{figure}

\smallskip
The global results in the one-dimensional case, obtained in \cite{Corcho-Matheus}, are based on a good control of the
$L^2$-norm of the solution $v$, which provides global well-posedness in $L^2\times L^2$. Global well-posedness
below $L^2$-regularity is obtained via the \emph{I-method} introduced by Colliander, Keel, Staffilani,
Takaoka and Tao in \cite{Imethod}. On the other hand, the global existence result, for any data in the space $H^1(\R^2)\times L^2(\R^2)$, established in
\cite{Corcho-Oliveira-Silva}, is obtained by using a careful estimate of the pseudo-Hamiltonian (\ref{Energy-2}) combined with the mass conservation (\ref{Conservation Law}) and the
Gagliardo-Nirenberg inequality in two dimensions:
\begin{equation}\label{Gagliardo-Nirenberg}
\|f\|_{L^4(\R^2)}\le c_2 \|f\|^{1/2}_{L^2(\R^2)}\|\nabla f\|^{1/2}_{L^2(\R^2)}.
\end{equation}

\smallskip
More recently, in \cite{Carvajal}, the authors showed global well-posedness for (\ref{SD}) for any initial data in the space $H^1(\R^2)\times H^1(\R^2)$ and for small data
in $H^s(\R^2)\times L^2(\R^2)$, with $2/3<s<1$, extending the previous results obtained in \cite{Corcho-Oliveira-Silva}.

\subsection{Main results}
We now present the new results obtained in this paper, for the Cauchy problem (\ref{SD}), in space dimensions $n=3$ and $n=4$.

Our first result is an addition to the local well-posedness results for $n=3$ established in \cite{Corcho-Oliveira-Silva} and concerns global well-posedness for the focusing case of (\ref{SD}),  with initial data in a broad subset of the space $H^1(\R^3)\times L^2(\R^3)$.

\smallskip
\begin{theorem}[\textbf{Global well-posedness in three  dimensions}]\label{GWP-Small-3D} Consider the system (\ref{SD}) with  $\lambda =-1$ and  initial data $(u_0,v_0)\in H^1({\R}^3)\times L^2({\R}^3)$, such that
the initial pseudo-energy $E_0:=E(0)$, given by \eqref{Energy-2}, is non-negative. Then, there exists a constant $\beta>0$, independent of the initial data,  
such that, if the initial data satisfies the condition\emph{:}
\begin{equation}\label{GWP-Small-3D-Small-Condittions-a}
	\|u_0\|^2_{L^2}E_0<\beta,
\end{equation}
then, a number $\gamma_0 \geq E_0$, depending on $\|u_0\|_{L^2}$ and  $E_0$, can be determined for which, if
\begin{equation}\label{GWP-Small-3D-Small-Condittions-b}
    \|\nabla u_0\|^2_{L^2}\leq \gamma_0,
\end{equation}
the local solution given by Theorem \ref{LWP-BourgainSpaces} can be extended to any time interval $[0,T]$.
\end{theorem}

\smallskip
\begin{remark} 
Regarding the previous theorem, we make the following three important observations.
 \begin{enumerate}[(a)]
\smallskip
\item If the initial pseudo-energy is negative, $E_0<0$, then the global control
of the $H^1({\R}^3)\times L^2({\R}^3)$ norm of the solution, as in the previous theorem, cannot be achieved. In
particular, as can be seen from \eqref{proof-gwp-small-3D-2} in the proof ahead, negative energy is incompatible
with the smallness condition \eqref{GWP-Small-3D-Small-Condittions-b} and in fact the unboundedness 
of $\|\nabla u\|_{L^2}$ is
not excluded in this case. In other words, values of $\|\nabla u_0\|_{L^2}$ in a small 
neighborhood of zero necessarily imply non-negative initial pseudo-energy $E_0$. 

\smallskip 
\item If $\lambda =-1$ and $v_0\geq0$ we have, from \eqref{Energy-2}, that $\|\nabla u_0\|^2_{L^2}\leq E_0\leq \gamma_0$, implying that, when condition \eqref{GWP-Small-3D-Small-Condittions-a} is satisfied, then \eqref{GWP-Small-3D-Small-Condittions-b} always is. Also note that the assumption \eqref{GWP-Small-3D-Small-Condittions-a} is not exactly a smallness condition on the data; it corresponds, rather, to a hyperbolic compensation between the energy $E_0$ and $\|u_0\|_{L^2}^2$, in which one of them can actually be large as long as the other is sufficiently small, so that the product satisfies
\eqref{GWP-Small-3D-Small-Condittions-a}. In fact, more generally, no matter what the sign of $v$ is, or
how large is its initial $L^2$ norm $\|v_0\|_{L^2}$, conditions \eqref{GWP-Small-3D-Small-Condittions-a}
and \eqref{GWP-Small-3D-Small-Condittions-b} can always be fulfilled by choosing 
$\|u_0\|_{H^1}$ sufficiently small.

\smallskip 
\item Recall that, for the defocusing $(\lambda =1)$ case of the cubic NLS (\ref{CNLS}), global well-posedness for any data in $H^1({\R}^3)$ is obtained using the fact that the conserved Hamiltonian
\begin{equation}\label{Hamiltonian-NLS}
\mathcal{H}(t)=\int_{\R^n}\Bigl(|\nabla u|^{2}+\lambda |u|^{4}\Bigl)dx=\mathcal{H}(0)
\end{equation}
is a positive quantity and thus provides an a priori estimate for the $H^1$-norm. Unfortunately, unlike in 
the cubic NLS case, neither is the pseudo-Hamiltonian (\ref{Energy-2}) of the Schr\"odinger-Debye system 
a conserved quantity, nor is it a positive quantity for either value of $\lambda$. So, 
besides the fact that, to the best of our knowledge,
no good control for this quantity is known, it is also not 
clear whether it would actually be helpful at all. Thus the  problem of global well-posedness for system (\ref{SD}), for arbitrary data in $H^1(\R^3)\times L^2(\R^3)$,
remains open in both cases $\lambda =\pm 1$. 
\end{enumerate}
\end{remark}

\smallskip
The critical Sobolev index for scaling invariance of the cubic NLS equation  (\ref{CNLS}) is given, as a function of the spatial dimension, by
\begin{equation}
s_n=\frac n2-1,
\end{equation}
from which it follows that $H^1$ is the critical Sobolev space in dimension $n=4$. 
Our remaining results all concern the Cauchy problem (\ref{SD}), precisely in four spatial dimensions and initial data 
in $H^1\times H^1$. 

\smallskip
We start by recalling the mixed $L^p$ norm notation, where $\|f\|_{L^p_IL_x^q}$ denotes the space-time norm  $$\|f\|_{L^p_IL_x^q}=\left(\int_I\|f(\cdot,t)\|^p_{L^q_x(\R^n)}dt\right)^{1/p},$$
for $I \subset \R_t$, some time interval.

\smallskip
\begin{theorem}[\textbf{Local well-posedness in  four dimensions}]\label{LWP-4D}
Given $(u_0, v_0) \in H^1({\R}^4)\times H^1({\R}^4)$, there exist
positive times $T_{\pm}=T_{\pm}(\mu, \|u_0\|_{H^1}, \|v_0\|_{H^1})$ and a unique solution to
the initial value problem (\ref{SD})
in the time interval $[-T_-, T_+]$ satisfying
\begin{equation}\label{LWP-4D-A}
(u, v)\in C\bigl([-T_-, T_+];\,H^1(\R^4)\times H^1(\R^4)\bigl),
\end{equation}
\begin{equation}\label{LWP-4D-B}
\|u\|_{L^{\infty}_IH_x^1}+ \|\nabla u\|_{L_I^2L_x^4}
+ \|v\|_{L^{\infty}_IH_x^1}<\infty.
\end{equation}
Moreover, for all\; $0<T'_{\pm}< T_{\pm}$, there exists a neighborhood $U'\times V'$ of $(u_0,v_0)$ in
$H^1(\R^4)\times H^1(\R^4)$ such that the map $(u_0,v_0) \longmapsto (u(\cdot, t), v(\cdot, t))$\; from $U'\times V'$ into the class defined by \eqref{LWP-4D-A}--\eqref{LWP-4D-B},
with $T'_{\pm}$ instead of $T_{\pm}$, is Lipschitz.
\end{theorem}

\smallskip
\begin{remark}\label{criticality}
The proof of Theorem \ref{LWP-4D} follows, without major difficulties, by adapting the standard techniques used to prove similar results for equation (\ref{CNLS}).
Notice, however, that whereas for the corresponding cubic NLS (\ref{CNLS}), in dimension four, the time of existence depends on the specific form of $u_0$ itself - a fact which is a typical feature of a truly critical problem (see \cite{Cazenave-Critical-1}) -
 we obtain here a local existence  result whose time of existence depends only on the size of the initial data $(u_0, v_0)$, that is, only on $\|u_0\|_{H^1}$  and  $\|v_0\|_{H^1}$. This can be interpreted as a regularizing effect
 introduced by the Debye delay equation, in \eqref{SD}, when compared to the $H^1$-critical cubic NLS \eqref{CNLS} for 
 $n=4$.
\end{remark}

\smallskip
Before stating the main results of this work, we point out a few important facts. First, we recall that the  existence of local solutions  $(u, v)$ for the Cauchy problem  \eqref{SD}, in the space $C\bigl((-T_-,\; T_+);\; H^s({\R}^n)\times H^s({\R}^n)\bigr)$, was established by B. Bid\'egaray in \cite{Bidegaray2}, for Sobolev indices $s>n/2$, 
using the algebra structure of the Sobolev spaces above that regularity index. Actually, the method of proof used in Theorem \ref{LWP-4D} can be similarly applied,
in four space dimensions $n=4$, to obtain local well-posedness for \eqref{SD} with initial data $(u_0, v_0)$ in  the space $H^s(\R^4)\times H^s(\R^4)$, for any integer Sobolev index $s \geq 1$. With that in mind, in the remaining part of this paper we will denote by $(-T_{l_s}, T_{r_s})$ the maximal time interval of existence of the corresponding 
solution $(u, v)$.

\smallskip 
Another important observation concerns the persistence property for system \eqref{SD}, that will be very useful in what follows. We make a precise statement of this property in the next remark.

\smallskip
\begin{remark}\label{persistence}
Let $(u_0, v_0) \in H^1(\R^4) \times H^1(\R^4)$ and $(u,v)$ be the corresponding solution given by Theorem \ref{LWP-4D}, defined 
in  $C\bigl((-T_{l_1}, T_{r_1});\,H^1(\R^4)\times H^1(\R^4)\bigl)$. Assuming furthermore that the initial 
data actually has higher regularity $(u_0, v_0) \in H^s(\R^4) \times H^s(\R^4)$, for some $s>1$, then the solution $(u,v)$ is also defined in the class $C\bigl((-T_{l_1}, T_{r_1});\,H^s(\R^4)\times H^s(\R^4)\bigl)$, i.e., $T_{r_s}= T_{r_1}$  and $T_{l_s}= T_{l_1}$.
\end{remark}
This phenomenon  of persistence of higher regularity is well known for NLS-type equations (see \cite{Linares-Ponce} pp.104, for example), where, by induction, one can proceed from $H^1$ to higher regularities, by
showing that the time of existence of the solution in $H^{k+1}$ is, at least, the same as for $H^{k}$, as long
as local well-posedness results are available for each such Sobolev index. It should be pointed out, though, that 
the step from $L^2(\R^4)=H^0(\R^4)$ to $H^1(\R^4)$ does not seem to hold for system \eqref{SD}, i.e. the time of existence of the $H^1(\R^4)$ solution cannot be proved to the same as the one in $L^2(\R^4)$ by following the same method.
 The reason for this discrepancy of persistence of regularity at the $L^2(\R^4)$ level is that we do not have at our disposal the required local well-posedness result, since if we were to replicate the proof of Theorem \ref{LWP-4D} for initial data in 
 $L^2(\R^4)\times L^2(\R^4)$ we would lose the time factor from the contraction scheme, yielding only a typical 
 existence result for small data. In a way, this somehow is more resembling of the $L^2$-critical behavior of the 
 cubic NLS in two dimensions rather than the related $H^1$-critical case in four spatial dimensions, that we are considering for system 
 \eqref{SD}, once again illustrating the aforementioned regularization effect introduced by the Debye relaxation.

\smallskip
The central ingredient used in the proof of the main result of this work is the following perturbed virial type identity.

\smallskip
\begin{theorem}[\textbf{Virial type identity}]\label{virial-theorem} Let $(u_0,v_0) \in H^s({\R}^n)\times H^s({\R}^n)$, with integer $s>1$ large enough (for $n=4$ it suffices to take $s=2$), and consider the corresponding $H^s\times H^s$-solution $(u, v)\in C\bigl((-T_{l_s}, T_{r_s});\, H^s(\R^n)\times  H^s(\R^n)\bigl)$ of (\ref{SD}) with $\lambda =-1$, defined on its maximal time interval $(-T_{l_s}, T_{r_s})$. Assume in addition that the initial variance is finite,
\begin{equation}\label{variance}
\int_{\R^n}|x|^2|u_0(x)|^2dx <\infty.
\end{equation}
Then, the function $t \longmapsto |\cdot|u(t,\cdot)$ is in $C\bigl((-T_{l_s}, T_{r_s});\, L^2(\R^n)\bigr)$, the function $t \longmapsto \displaystyle \int_{\R^n}|x|^2|u(x,t)|^2dx$ is in $C^2(-T_{l_s}, T_{r_s})$ and we have
\begin{equation}\label{virial-a}
\frac{d}{dt}\frac12\int_{\R^n}|x|^2|u|^2dx=\emph{Im}\int_{\R^n}(x\cdot \nabla u)\bar{u}\,dx
\end{equation}	and
\begin{equation}\label{virial-b}
\frac{d^2}{dt^2}\frac12 \int_{\R^n}|x|^2|u|^2dx=E(t)+(n-2)\int_{\R^n}v|u|^2\,dx-\int_{\R^n}v^2\,dx
+\int_{\R^n}( x\cdot \nabla |u|^2) v\,dx.
\end{equation}
\end{theorem}

\smallskip
We now state the principal theorem in this work which establishes that, for the focusing case of system (\ref{SD}) in dimension $n=4$ and fixed radial and smooth initial data, there 
is a time $T_0$, that only depends on the data, for which there cannot exist solutions on $[0,T]$,
with $T>T_0$,
which have a higher regularity Sobolev norm that is uniformly bounded in the delay parameter $\mu$. This turns into an alternative result, where either the existence, or the uniform boundedness, fails. 
The main requirement for the initial data, similarly to the usual blow-up condition for $L^2$-critical or supercritical  NLS equations, is that the energy be negative, which is something that can be 
achieved within the Schwartz class, as will be seen afterwards. The theorem reads as follows: 

\smallskip
\begin{theorem}[\textbf{Blow-up or unboundedness of Sobolev norms for small $\mu$}]\label{blowup-theorem}
Consider radial initial data $(u_0,v_0) \in \mathscr{S}({\R}^4)\times \mathscr{S}({\R}^4)$ and let 
$(u_{\mu}, v_{\mu})\in C\bigl((-T_{{\mu}_{l_1}}, T_{{\mu}_{r_1}});\, H^1(\R^4)\times  H^1(\R^4)\bigl)$ be the corresponding maximal time interval $H^1\times H^1$-solution  of \eqref{SD} with $\lambda=-1$. Assume, in addition, that the initial data $(u_0, v_0)$ is such that 
its pseudo-energy satisfies
\begin{equation}\label{blowup-cond}
E_0:=E(0) <0.
\end{equation}
Then, there exists a positive time  $T_0$, depending only on the data, such that 
one of the following two alternatives holds. Either
\begin{enumerate}
\item[\textbf{(a)}] there is a vanishing sequence  of values of the delay parameter $\mu\to 0$ for which the corresponding solutions all
satisfy $T_{{\mu}_{r_1}} \leq T_0<\infty$, that is, these $H^1\times H^1$-solutions blow-up at a finite time
in the interval $(0, T_0]$  
\end{enumerate}
or 
\begin{enumerate}
\item[\textbf{(b)}] there is $T\geq T_0$ and $0<\mu_0<1$ for which the initial value problem \eqref{SD} has solutions in $[0,T]$, for all $0<\mu\leq\mu_0$,  $(u_{\mu}, v_{\mu})\in C\bigl((-T_{{\mu}_{l_1}}, T_{{\mu}_{r_1}});\, H^1(\R^4)\times  H^1(\R^4)\bigl)$, with $T_0\leq T<T_{{\mu}_{r_1}}$, in which case, 
for all $s> 4 $ the Sobolev norms of the solutions for the Schr\"odinger component of the system are unbounded when $\mu \to 0$, that is, 
$$\sup\limits_{0< \mu \le \mu_0} \Big(\sup\limits_{0\le t \le T}\|u_{\mu}(\cdot, t)\|_{H^s}\Big)=+\infty.$$
\end{enumerate}
In case blow-up does occur, which can also happen for alternative (b) at times larger than $T_0$, then 
$$
\lim\limits_{t\nearrow T_{{\mu}_{r_1}}}\|(u_{\mu}(\cdot, t),v_{\mu}(\cdot, t))\|_{H^1(\R^4)\times  H^1(\R^4)}=
+\infty,
$$
and in particular this implies that $\lim\limits_{t\nearrow T_{{\mu}_{r_1}}}\|\nabla v_{\mu} (\cdot,t)\|_{L^2(\R^4)}= +\infty$  and also
$$\lim\limits_{t\nearrow T_{{\mu}_{r_1}}}\|u_{\mu}(\cdot,t)\|_{L^{\infty}(\R^4)}= \lim\limits_{t\nearrow T_{{\mu}_{r_1}}}\|\nabla u_{\mu} (\cdot,t)\|_{L^2(\R^4)}=+\infty.$$
\end{theorem}

\smallskip
\begin{remark}
We note the following three  observations.
\begin{enumerate}[(a)]
\smallskip
\item For fixed initial data,  alternative (a) in the previous theorem only guarantees blow-up of the solutions for a sequence of small enough parameters $\mu$, converging to zero. Of course, for any other values of $\mu$ blow-up solutions can then be obtained from these by using the rescaling \eqref{parameter-rescaling}, but it should not be forgotten that this also rescales the initial data, changing it accordingly.

\smallskip
\item Our proof of the virial identity exploits ideas similar to the ones developed by F. Merle \cite{Merle}, 
in the context of the Zakharov system. Also, the method of proof of Theorem \ref{blowup-theorem} remains valid for any dimension $n\geq 4$ as long as a local well-posedness result in $H^1 \times H^1$ is available for such a dimension. In
that case, the unboundedness of Sobolev norms in alternative (b) would hold for regularities $s>2+n/2$.

\smallskip
\item The time $T_0$ seems to signal the onset of some kind of pathological behaviour of the solutions to \eqref{SD}, as $\mu \to 0$. The
	alternative scenario of Theorem \ref{blowup-theorem} leaves open several possibilities. It might be the case
	that the solutions have uniformly bounded Sobolev norms on compact intervals $[0,T]$, for any $T<T_0$,
	and that they all blow-up exactly at $T_0$. Or, it might just happen that there really is no blow-up in finite time and that
	all solutions are global for positive times, except that the higher regularity Sobolev norms become unbounded 
	in $\mu$, as $\mu \to 0$, once the time $T_0$ is crossed. In any case, all of this might also be suggestive
	of limiting pathological behavior of the solutions, as $\mu \to 0$, hinting that the solution 
	to the cubic NLS in four dimensions, for the same initial data, could also exhibit some type of singular properties
	at $T_0$.
\end{enumerate}
\end{remark}

\smallskip
For high regularity initial data, with small $H^1 \times H^1$ norm we can, nevertheless, establish the following global well-posedness result.

\smallskip
\begin{proposition}[\textbf{Global well-posedness in four dimensions}]\label{GWP-Small-4D}
Consider the system (\ref{SD}) with  $\lambda =-1$ and  initial data $(u_0,v_0)\in H^s({\R}^4)\times H^s({\R}^4)$, for some $s>2$, such that
the initial pseudo-energy $E_0:=E(0)$, given by \eqref{Energy-2}, is non-negative. Then, there exists a constant $\beta>0$, independent of the initial data,  
such that, if the initial pseudo-energy satisfies\emph{:}
\begin{equation}
E_0<\beta,
\end{equation}
then, a number $\gamma_0 \geq 0$, depending only on  $E_0$, exists for which, if
\begin{equation}
\|\nabla u_0\|_{L^2}\leq \sqrt{\gamma_0},
\end{equation}
the local solution given by Theorem \ref{LWP-BourgainSpaces} can be extended to any time interval $[0,T]$.
\end{proposition}

\smallskip
We finish this section by establishing the existence of  functions that satisfy the hypotheses for the initial data in Theorem \ref{blowup-theorem} 
providing  explicit examples for which they are applicable.

\smallskip
\begin{proposition}
There exist functions $(u_0,v_0) \in \mathscr{S}({\R}^4)\times \mathscr{S}({\R}^4)$ such that
\begin{equation}\label{energianegativa}
\int_{\R^4}\Bigl(|\nabla u_0|^{2}+2v_0|u_0|^2 + v_0^2\Bigl)dx<0.
\end{equation}
\end{proposition}

\noindent \textbf{Proof.} Take $\phi \in C^{\infty}_c(\R)$ such that
$$
\phi(s)=
\begin{cases}
1 \quad \mbox{if}\quad s\leq 1,\\
0 \quad \mbox{if}\quad s \geq 2,
\end{cases}
$$
and make 
$$
u_0=\frac1N\, \phi\left( \frac{|x|}{N^2}\right) \qquad \mbox{and} \qquad v_0=-|u_0|^2,
$$
with large $N$ to be chosen conveniently at the end. Of course, $u_0, v_0 \in \mathscr{S}({\R}^4)$.

Computing the gradient of $u_0$ we obtain
$$
\nabla u_0(x)=\frac{1}{N^3}\, \phi'\left( \frac{|x|}{N^2}\right) \frac{x}{|x|},
$$
noting that the apparent singularity at $x=0$ does not pose any problem as $\phi=1$ in a neighborhood of the origin and therefore
$\phi'=0$ in that same neighborhood.

Gathering everything in the integral formula \eqref{energianegativa} we obtain
$$
\int_{\R^4}\Bigl(|\nabla u_0|^{2}+2v_0|u_0|^2 + v_0^2\Bigl)dx=\int_{\R^4} \frac{1}{N^6} \left| \phi'\left( \frac{|x|}{N^2}\right)  \right| ^2 - \frac{1}{N^4}\left| \phi\left( \frac{|x|}{N^2}\right)  \right| ^4 dx,
$$
which, using spherical coordinates, can be written as
$$
\omega_3 \int_0^{\infty}\rho^3 \left( \frac{1}{N^6} \left| \phi'\left( \frac{\rho}{N^2}\right)  \right| ^2 - \frac{1}{N^4}\left| \phi\left( \frac{\rho}{N^2}\right)  \right| ^4 \right) d\rho,
$$
where $\omega_3$ is the area of the unit three dimensional sphere $\mathbb{S}^3 \subset \R^4$. Now, doing a
change of variables $z=\frac{\rho}{N^2}$, we finally obtain
$$
\omega_3 \int_0^{\infty}N^6 z^3 \left( \frac{1}{N^6} \left| \phi'(z)  \right| ^2 - \frac{1}{N^4}\left| \phi(z)  \right| ^4 \right) N^2 dz= N^2 \omega_3  \int_0^{\infty} z^3 \left| \phi'(z)  \right| ^2 dz -N^4 \omega_3  \int_0^{\infty} z^3 \left| \phi(z)  \right| ^4 dz,
$$
from which we conclude that, by choosing $N$ large enough, this quantity can be made
negative. $\square$

\subsection{Cubic NLS versus Schr\"odinger-Debye}
As pointed out in \cite{Corcho-Oliveira-Silva}, as well as in the Remark \ref{criticality} above, concerning the
criticality of the local well-posedness result for $n=4$, the delay term $\mu v_t$ in (\ref{SD}) induces a regularization
with respect to the flow of the corresponding cubic NLS. We summarize, in the following table, a comparison of the known results concerning the local well-posedness for these equations.

\smallskip
\begin{table}[h]\caption{Local well-posedness ($\lambda =\pm 1$)}
\begin{tabular}{|c|l|l|}
\rowcolor{yellow!25}
\hline
\parbox[c][0.04\linewidth][c]{0.04\linewidth}{\centering $n$} & \bf Cubic NLS in $\boldsymbol{H^s(\R^n)}$ & \bf Schr\"odinger-Debye in $\boldsymbol{H^s(\R^n)\times H^{\kappa}(\R^n)}$\\
\hline
\hline
\parbox[c][0.04\linewidth][c]{0.04\linewidth}{\centering $1$} & $s\ge 0$ (\cite{Cazenave-Critical-2, Ginibrevelo,Y-Tsutsumi})
& $|s|-\frac 12\le \kappa <\min\bigl\{s+\frac 12,\; 2s+\frac 12\bigl\}$ (\cite{Corcho-Matheus})\\
\hline
\parbox[c][0.04\linewidth][c]{0.04\linewidth}{\centering $2$} &  $s\ge 0$ (\cite{Cazenave-Critical-1, Cazenave-Critical-2, Ginibrevelo}) &
    $\max \bigl\{0, s-1\bigl\} \le \kappa \le \min\{2s,\; s+1\}$ (\cite{Corcho-Oliveira-Silva})\\
\hline
\parbox[c][0.04\linewidth][c]{0.04\linewidth}{\centering $3$} &   $s\ge \frac 12$ (\cite{ Cazenave-Critical-2, Ginibrevelo})&  $\max \{0, s-1\} \le \kappa \le \min\{2s,\; s+1\}$ (\cite{Corcho-Oliveira-Silva})\\
\hline
\parbox[c][0.04\linewidth][c]{0.04\linewidth}{\centering $4$} & $s \geq 1$ (\cite{ Cazenave-Critical-1, Cazenave-Critical-2, Ginibrevelo}) & $(s,\kappa)=(1,1)$\\
\hline
\end{tabular}
\end{table}

\bigskip
The plan of the paper is the following. In Section \ref{Section-sgwp-3D}, we prove global well-posedness in the space $H^1(\R^3)\times L^2(\R^3)$, under certain restrictions for the initial data. In Section \ref{Section-wp-4D}, we establish the local theory in $H^1(\R^4)\times H^1(\R^4)$. In Section \ref{Section-virial-identity} we derive the virial type  identity \eqref{virial-b} from which, in Section \ref{Section-blowup} the proof of Theorem \ref{blowup-theorem} follows, by using a contradiction argument. Finally, in Section \ref{GWP-4D} we prove
global-wellposedness in $H^1(\R^4)\times H^1(\R^4)$, for higher  regularity initial data and small $H^1$ norm.

% % % % % % % % % % % % % % % % % % % % % % % % % % % % % % % % % % % % % % % % % % % % % % % % % % % % % % % % % % % % % % % % % % % %
\section{\bf{Global Well-posedness in Three Dimensions.}}\label{Section-sgwp-3D} 
In this section we derive a priori estimates in the spaces  $H^1(\R^3)\times L^2(\R^3)$ for the focusing case of (\ref{SD}), from which the global well-posedness follows for initial data satisfying conditions
\eqref{GWP-Small-3D-Small-Condittions-a} and \eqref{GWP-Small-3D-Small-Condittions-b}. The following version of the Gagliardo-Nirenberg inequality, for $n=3$ will be used in the proof
\begin{equation*}
\|f\|_{L^4(\R^3)}\le c_3 \|f\|^{1/4}_{L^2(\R^3)}\|\nabla f\|^{3/4}_{L^2(\R^3)}.
\end{equation*}

\smallskip
\subsection{Proof of Theorem \ref{GWP-Small-3D}}
 
Let $n=3$, $\lambda=-1$ and consider the $H^1(\R^3)\times L^2(\R^3)$ solution  $\bigl(u(\cdot, t),v(\cdot,t)\bigl)$ of (\ref{SD}) established in Theorem \ref{LWP-BourgainSpaces} and defined on its maximal positive time
interval $\bigl[0, T_*\bigr).$ Using the pseudo-energy (\ref{Energy-2}) and combining the H\"older, Gagliardo-Nirenberg and Young inequalities we have
\begin{equation}
\begin{split}\label{proof-gwp-small-3D-1}
\|\nabla u(\cdot, t)\|_{L^2}^2 + \|v(\cdot, t)\|_{L^2}^2&\le E_0 -2\int_{\R^3}v(\cdot, t)|u(\cdot, t)|^2dx\\
&\le E_0 + 2 \|v(\cdot, t)\|_{L^2}\|u(\cdot, t)\|^2_{L^4}\\
&\le E_0 + 2 c_3^2 \|v(\cdot, t)\|_{L^2}
\|u(\cdot, t)\|^{\frac{1}{2}}_{L^2}\|\nabla u(\cdot, t)\|^{\frac{3}{2}}_{L^2}\\
&\le E_0 + \|v(\cdot, t)\|^2_{L^2} + c_3^4\|u(\cdot, t)\|_{L^2}\|\nabla u(\cdot, t)\|^3_{L^2}.
\end{split}
\end{equation}

\smallskip
Now, we define the continuous function
$$\phi(t):=\|\nabla u(\cdot, t)\|_{L^2}^2,$$
for all $0\le t < T_*$. Using the conservation of the $L^2$-norm of the solution $u(\cdot, t)$,
the inequality (\ref{proof-gwp-small-3D-1}) yields the a priori estimate
\begin{equation}\label{proof-gwp-small-3D-2}
0\le \phi(t) \le E_0 +  \nu_0\phi(t)^{\frac 32}\quad \text{with}\quad \nu_0=c_3^4 \|u_0\|_{L^2}.
\end{equation} 

\smallskip
If $u_0=0$ then the solution is obviously global in time as, in that case, $u(x,t)=0$ and $v(x,t)=
e^{-t/\mu}v_0(x)$ for all $(x,t) \in \R^3 \times \R$. So we can assume, for the remaining part of the proof,
that $\|u_0\|_{L^2} >0$.

\smallskip
In Figure \ref{Fig-Bootstrap} below, we draw the graph of the convex function $f(x)= E_0+\nu_0x^{\frac32}$,
indicating its point $\bigl(x_0, f(x_0)\bigr)$ of slope one, at $x_0=\frac{4}{9\nu_0^2}$, with its tangent
line.

\smallskip
\begin{figure}[ht]
\begin{center}
\begin{tikzpicture}[domain=0:5]
%\draw[very thin,color=gray] (-0.1,-1.1) grid (6.4,3.9);
\draw[->] (-0.5,0)--(5.5,0) node[below] {$x$};
\draw[->] (0,-1)--(0,5.5) node[right] {$y$};
\draw[thick] plot (\x, \x)  node at (2.5,3) {{\small $y=x$}}; 
\draw[thick, domain=0:5] plot (\x, 1 + 0.18*\x^2) node at (4.3,5) {{\small $f(x)$}};
\draw[thin, color=red] (0,-0.43)--(5,4.57); %node at (4.8,4) {{\small $L$}};
\draw[thick, dashed] (2.9,0)--(2.9,2.5)--(0,2.5) node at (2.9,2.5){{\small $\bullet$}};
\draw[thick, dashed] (1.32,0)--(1.32,1.32)--(0,1.32) node at (1.32,1.32){{\small $\bullet$}};
\draw[thick, dashed] (4.25,0)--(4.25,4.25)--(0,4.25) node at (4.25,4.25){{\small $\bullet$}};
\node at (0,1){{\small $\bullet$}};
\node at (-0.4,1){{\small $E_0$}};
\node at (1.32,0){{\small $\bullet$}};
\node at (1.32,-0.3){{\small $\gamma_0$}};
\node at (4.25,0){{\small $\bullet$}};
\node at (4.25,-0.3){{\small $\tilde{\gamma}_0$}};
\node at (2.9,0){{\small $\bullet$}};
\node at (2.9,-0.3){{\small $x_0$}};
%\node at (0,-0.43){{\small $\bullet$}};
%\node at (-0.3,-0.43){{\small $\eta_0$}};
\end{tikzpicture}
\end{center}
\vspace{-0.3cm}
\caption{}
\label{Fig-Bootstrap}
\end{figure}
We observe that, if the initial data is such that the condition $f(x_0)<x_0 \Leftrightarrow E_0< \frac{4}{27\nu_0^2}$ is satisfied, which is equivalent to 
\begin{equation}\label{cond-eta0}
\|u_0\|^2_{L^2}\, E_0< \frac{4}{27c_3^8},
\end{equation}
then, the function $f$ intersects the line $y=x$ at two points, $x=\gamma_0$ and $x=\tilde{\gamma}_0$, with $\gamma_0<x_0<\tilde{\gamma}_0$. Thus, if 
\begin{equation}\label{small-data-cond}
\phi(0)= \|\nabla u_0\|^2_{L^2}\leq \gamma_0
\end{equation}
then $0\le \phi(t) \le f(\phi(t))$ and using the continuity of the function $\phi: \bigl[0, T_*\bigr) \longrightarrow \R$ we have that the values of the function $\phi$ are trapped in the interval
\begin{equation}\label{proof-gwp-small-3D-4}
0\le  \phi(t) \le  \gamma_0, \quad \text{for all}\quad 0\le t < T_*,
\end{equation}
according to Figure \ref{Fig-Bootstrap}.

\smallskip
As for the term $\|v(\cdot, t)\|_{L^2}$, we return to \eqref{proof-gwp-small-3D-1} and repeat the last step, using Young's inequality with any small $0<\varepsilon<1$,
$$\|\nabla u(\cdot, t)\|_{L^2}^2 + \|v(\cdot, t)\|_{L^2}^2 \le
E_0 + \varepsilon\|v(\cdot, t)\|^2_{L^2} + \frac{c_3^4}{\varepsilon}\|u(\cdot, t)\|_{L^2}\|\nabla u(\cdot, t)\|^3_{L^2},$$
which, by using the conservation of the $L^2$-norm of the solution $u$, implies 
$$\|v(\cdot, t)\|_{L^2}^2 \leq \frac{E_0}{1-\varepsilon}+\frac{c_3^4}{\varepsilon(1-\varepsilon)}\|u_0\|_{L^2}\|\nabla u(\cdot, t)\|^3_{L^2} \leq 
\frac{E_0}{1-\varepsilon}+\frac{c_3^4}{\varepsilon(1-\varepsilon)}\|u_0\|_{L^2} \gamma_0^{\frac32},$$
for all $0\le t < T_*$.

\smallskip
We thus conclude that $T_*$ must be infinite, since from local theory we know that, if $T_*< \infty$, then $\|\nabla u(\cdot, t)\|_{L^2}^2 + \|v(\cdot, t)\|_{L^2}^2$
would necessarily have to blow-up at this endpoint, and we have proved that both terms remain bounded.
Finally, combining  \eqref{cond-eta0}  and \eqref{small-data-cond}, and defining $\beta=\frac{4}{27c_3^8}$,
 we obtain the conditions 
\eqref{GWP-Small-3D-Small-Condittions-a}-\eqref{GWP-Small-3D-Small-Condittions-b}. The proof is thus finished. $\square$

% % % % % % % % % % % % % % % % % % % % % % % % % % % % % % % % % % % % % % % % % % % % % % % % % % % % % % % % % %
\section{\bf{Local Theory in} $\boldsymbol{H^1(\mathbb{R}^4)\times H^1(\mathbb{R}^4)}$.}\label{Section-wp-4D}

In this section we present the proof of Theorem \ref{LWP-4D}. We recall  the Strichartz estimate 
for the free Schr\"odinger group $\displaystyle S(t)=e^{it\Delta/2}$ in the euclidean space $\R^4$.

\smallskip
\begin{lemma}[\textbf{Strichartz estimates} \cite{Cazenave-Book, Keel-Tao}] Let $(p_1,q_1)$ and $(p_2,q_2)$ be two pairs of admissible exponents for $S(t)=e^{it\Delta/2}$ in $\R^4$; that is, both satisfying
	the condition
	\begin{equation}\label{Strichart-Admissible}
	\frac{2}{p_i}=4\left(\frac12 -\frac 1{q_i} \right) \quad \text{and} \quad 2\le q_i\le 4\quad (i=1,2).
	\end{equation}
	Then, for any $0<T\leq \infty$, we have
	\begin{equation}\label{Strichart-homogeneous}
	\|S(t)f\|_{L^{p_1}_TL^{q_1}_x}\le c\|f\|_{L^2(\R^4)},
	\end{equation}
	as well as the non-homogeneous version
	\begin{equation}\label{Strichart-non-homogeneous}
	\left\|\int_0^tS(t-t')g(\cdot, t')dt'\right\|_{L^{p_1}_TL_x^{q_1}}\le c \|g\|_{L^{p'_2}_TL_x^{q'_2}},
	\end{equation}
	where $1/p_2+1/p_2'=1$, $1/q_2 + 1/q'_2=1$ and $\|f\|_{L^p_TL_x^q}=\|f\|_{L^p_{[0,T]}L_x^q}$. The constants
	in both inequalities are independent of $T$. 
\end{lemma}

\smallskip
\subsection{Proof of Theorem \ref{LWP-4D}}
Consider de integral formulation  for (\ref{SD}), given by
$$
\left\lbrace 
\begin{array}{rcl}
u(\cdot, t)&\!\!\!\! = \!\!\!\!& \displaystyle S(t)u_0-i\int_0^tS(t-\tau)u(\cdot, \tau)v(\cdot, \tau)d\tau,\\[0.3cm]
v(\cdot, t)&\!\!\!\! =\!\!\!\! & \displaystyle e^{-t/\mu}v_0+\tfrac{\lambda}{\mu}\int_0^t e^{-\frac{t-\tau}{\mu}}|u(\cdot, \tau)|^2d\tau,
\end{array}
\right.
$$
from which we define the two operators
\begin{align}
&\Phi_1(u,v):=S(t)u_0-i\int_0^tS(t-\tau)u(\cdot, \tau)v(\cdot, \tau)d\tau,\label{Sd-Integral-u}\\
&\Phi_2(u,v):=e^{-t/\mu}v_0+\tfrac{\lambda}{\mu}\int_0^t e^{-\frac{t-\tau}{\mu}}|u(\cdot, \tau)|^2d\tau,\label{Sd-Integral-v}
\end{align}
and the sets
\begin{align}
&U_{\rho_1,T}=\Bigl\{u: [0,T]\times \R^4 \rightarrow \C;\; \|u\|_U:=\|u\|_{L_T^{\infty}H^1_x} + \|\nabla u\|_{L^2_TL^4_x}\le \rho_1\Bigl\}\\
\intertext{and}
&V_{\rho_2,T}=\Bigl\{v: [0,T]\times \R^4 \rightarrow \R;\; \|v\|_V:=\|v\|_{L_T^{\infty}H^1_x}\le \rho_2\Bigl\}.
\end{align}

As usual, we will next choose $\rho_1$, $\rho_2$ and $T$ so that the operator $\Phi=(\Phi_1, \Phi_2)$ maps
$U_{\rho_1,T}\times V_{\rho_2,T}$ to itself, 
$$\Phi=(\Phi_1, \Phi_2):U_{\rho_1,T}\times V_{\rho_2,T} \longrightarrow  U_{\rho_1,T}\times V_{\rho_2,T},$$
and is a contraction, with the norm
\begin{equation}\label{proof-lwp-4d}
\|(u,v)\|_{U\times V}=\|u\|_U + \|v\|_V,
\end{equation}
yielding the fixed point that satisfies 
the integral formulation of the problem. 

Indeed, note that
\begin{equation*}\label{proof-lwp-4d-a}
\begin{split}
\|\Phi_1(u,v)\|_U&\le \|S(t)u_0\|_U+ \Bigl \|\int_0^tS(t-\tau)u(\cdot, \tau)v(\cdot, \tau)d\tau \Bigl\|_U\\
&\le c \|u_0\|_{H^1} + c\Bigl(\|uv\|_{L^2_TL^{4/3}_x} + \|\nabla(uv)\|_{L^2_TL^{4/3}_x}\Bigl).
\end{split}
\end{equation*}
This follows, for the homogeneous term, by (\ref{Strichart-homogeneous}) with $(p_1,q_1)=(\infty, 2)$ and $(p_1,q_1)=(2,4)$. For the non-homogeneous
term we used (\ref{Strichart-non-homogeneous}) with the same two pairs of $(p_1,q_1)$, chosen in the previous case, and with $(p'_2,q'_2)=(2,4/3)$. Now, using  H\"older and Sobolev inequalities we obtain, for all $(u,v)\in U_{\rho_1,T}\times V_{\rho_2,T}$, the following estimates:
\begin{equation}\label{proof-lwp-4d-b}
\begin{split}
\|\Phi_1(u,v)&\|_U\le c \|u_0\|_{H^1} +
c\|u\|_{L^2_TL^4_x}\left(\|v\|_{L^{\infty}_TL^2_x} + \|\nabla v\|_{L^{\infty}_TL^2_x}\right) + c\|v\|_{L^2_TL^4_x}\|\nabla u\|_{L^{\infty}_TL^2_x}\\
& \le c \|u_0\|_{H^1} + c\|\nabla u\|_{L^2_TL^2_x}\|v\|_{L^{\infty}_TH^1_x} + c\|\nabla v\|_{L^2_TL^2_x}\|\nabla u\|_{L^{\infty}_TL^2_x}\\
&\le c \|u_0\|_{H^1}+ c\sqrt{T}\|u\|_{L^{\infty}_TH^1_x}\|v\|_{L^{\infty}_TH^1_x}\\
&\le c \|u_0\|_{H^1}+ c\sqrt{T}\rho_1\rho_2.
\end{split}
\end{equation}
On the other hand, applying the Minkowski and H\"older inequalities to \eqref{Sd-Integral-v} we get
\begin{equation*}\label{proof-lwp-4d-c}
\|\Phi_2(u,v)\|_{H^1_x}\le e^{-t/\mu}\|v_0\|_{H^1} + \tfrac{1}{\mu}\int_0^t e^{-\frac{t-\tau}{\mu}}\left( \|u\bar{u}\|_{L^2_x}
+2\|\bar{u}\nabla u\|_{L^2_x} \right)d\tau.\\
\end{equation*}
Again, using H\"older and Sobolev inequalities, we have that
\begin{equation}\label{proof-lwp-4d-d}
\begin{split}
\|\Phi_2(u,v)&\|_{H^1_x}\le e^{-t/\mu}\|v_0\|_{H^1} + \tfrac{1}{\mu}\int_0^t e^{-\frac{t-\tau}{\mu}}\left( \|u\|^2_{L^4_x}
+2 \|u\|_{L^4_x}\|\nabla u\|_{L^4_x}\right)d\tau\\
&\le e^{-t/\mu}\|v_0\|_{H^1} + \tfrac{1}{\mu}\left(\int_0^t e^{-2\frac{t-\tau}{\mu}}d\tau\right)^{\frac12}
\left( \|u\|_{L^{\infty}_{[0,t]}L^4_x} \|u\|_{L^2_{[0,t]}L^4_x}
+ 2\|u\|_{L^{\infty}_{[0,t]}L^4_x} \|\nabla u\|_{L^2_{[0,t]}L^4_x}\right)\\
&\le e^{-t/\mu}\|v_0\|_{H^1} + \sqrt{\frac{1-e^{-\frac{2t}{\mu}}}{2\mu}}
\left( \sqrt{t}\|\nabla u\|^2_{L^{\infty}_{[0,t]}L^2_x}+ 2\|\nabla u\|_{L^{\infty}_{[0,t]}L^2_x}
\|\nabla u\|_{L^2_{[0,t]}L^4_x}\right)\\
&\le e^{-t/\mu}\|v_0\|_{H^1} + c\frac{\sqrt{t}}{\mu}
\left( \sqrt{t}\|\nabla u\|^2_{L^{\infty}_{[0,t]}L^2_x}+ 2\|\nabla u\|_{L^{\infty}_{[0,t]}L^2_x}
\|\nabla u\|_{L^2_{[0,t]}L^4_x}\right).
\end{split}
\end{equation}
Thus, for all $u\in U_{\rho_1,T}$, it follows that
\begin{equation}\label{proof-lwp-4d-e}
\|\Phi_2(u,v)\|_{L^{\infty}_TH^1_x}\le \|v_0\|_{H^1} + c\frac{\sqrt{T}}{\mu}(\sqrt{T} +2)\rho_1^2.
\end{equation}
Now, if we fix  $\rho_1=2c \|u_0\|_{H^1}$ and $\rho_2=2\|u_0\|_{H^1}$ and take $T>0$ such that
\begin{equation}\label{proof-lwp-4d-timecondition-1}
c\sqrt{T}\rho_2\le \frac12\quad \text{and} \quad c\frac{\sqrt{T}}{\mu}(\sqrt{T} +2)\rho_1^2\le \rho_2,
\end{equation}
it follows that the application $\Phi=(\Phi_1, \Phi_2)$ is well-defined and
$\Phi\bigl(U_{\rho_1,T}\times V_{\rho_2,T}\bigl)\subset U_{\rho_1,T}\times V_{\rho_2,T}$.

\smallskip
The same type of estimates used in (\ref{proof-lwp-4d-b}) and 
(\ref{proof-lwp-4d-d}) show that $\Phi=(\Phi_1, \Phi_2)$ is also a contraction in $U_{\rho_1,T}\times V_{\rho_2,T}$ (with, eventually, smaller choices for $T, \rho_1$ and $\rho_2$) and this concludes
the proof. $\square$

\smallskip
\begin{remark}\label{timereversal}
	The local theory  obtained in the previous proof for the time interval $[0,T]$ can be extended to an interval $[-\tilde{T},0]$ by reflection of
	the time variable. More precisely, we can establish a similar local theory  in $[0,\tilde{T}]$ for the auxiliary system
	\begin{equation}\label{SD-auxiliary}
	\begin{cases}
	i\tilde{u}_t-\frac{1}{2}\Delta \tilde{u} = \tilde{u}\tilde{v},\\
	\mu\tilde{v}_t - \tilde{v}=\lambda|\tilde{u}|^2,\\
	\tilde{u}(x,0)=u_0(x),\quad \tilde{v}(x,0)=-v_0(x),
	\end{cases}
	\end{equation}
	by just slightly modifying the estimates made in
	(\ref{proof-lwp-4d-d}). Indeed, while \eqref{proof-lwp-4d-b} follows exactly in the same way, the solution 
	for the ODE in $\tilde{v}$ now yields an
	integral term $\displaystyle \int_0^t e^{2\frac{t-\tau}{\mu}}d\tau$ in (\ref{proof-lwp-4d-d}) 
	that can be estimated as follows:
	$$\left(\int_0^t e^{2\frac{t-\tau}{\mu}}d\tau\right)^{\frac12}\le e^{t/\mu}\sqrt{t}\le e^{\tilde{T}/\mu}\sqrt{\tilde{T}},$$
	for all $t\in [0,\tilde{T}]$. Thus, the remaining estimates also follow in a similar way as before. Then,
	$$u(x,t):=\tilde{u}(x,-t)\quad \text{and}\quad v(x,t):=-\tilde{v}(x,-t)$$
	are local solutions for the Cauchy problem (\ref{SD}), with the same initial data $(u_0, v_0)$, in the time interval $[-\tilde{T},0]$.
\end{remark}

% % % % % % % % % % % % % % % % % % % % % % % % % % % % % % % % % % % % % % % % % % % % % % % % % % % % % % % %

\section{\bf{Virial Type Identity for the Schr\"odinger-Debye System.}}\label{Section-virial-identity}

The following result is the main ingredient in the proof of Theorem \ref{blowup-theorem},  
for solutions of the focusing case of \eqref{SD}, in four spatial dimensions.

\subsection{Proof of Theorem \ref{virial-theorem}}
First, we prove (\ref{virial-a}) by following a similar argument as in the analogous result for the NLS equation. Multiplying the first equation of (\ref{SD}) by $|x|^2\bar{u}$,
integrating in the $x$ variable and taking the imaginary part, we obtain
\begin{equation}\label{proof-virial-a0}
\begin{split}
\frac{d}{dt}\int_{\R^n}|x|^2|u|^2dx&=-\text{Im}\int_{\R^n}|x|^2\bar{u}\Delta u\,dx\\
&=\text{Im}\int_{\R^n}\Bigl(|x|^2|\nabla u|^2+ 2\bar{u}(x\cdot \nabla u)\Bigl)dx\\
&=2\, \text{Im}\int_{\R^n}(x\cdot \nabla u)\bar{u}dx,
\end{split}
\end{equation}
which yields (\ref{virial-a}). This formal procedure can be made rigorous through a regularizing technique (see \cite{Cazenave-Book}).

In order to prove (\ref{virial-b}) we need to compute the term $\displaystyle \text{Im}\frac{d}{dt}\int_{\R^n}(x\cdot \nabla u)\bar{u}\,dx$.
Again, proceeding formally assuming all computations can be performed, we start by rewriting the derivative in time as follows:
\begin{equation}\label{proof-virial-a}
\frac{d}{dt}\int_{\R^n}(x\cdot \nabla u)\bar{u}\,dx=\int_{\R^n}(x\cdot \nabla u)\bar{u}_t\,dx + \int_{\R^n}(x\cdot \nabla u_t)\bar{u}\,dx.
\end{equation}
Now, using integration by parts we get
\begin{equation}\label{proof-virial-b}
\begin{split}
\int_{\R^n}(x\cdot \nabla u_t)\bar{u}\,dx
&=-n\int_{\R^n}u_t\bar{u}\,dx -\int_{\R^n}(x\cdot \nabla \bar{u})u_t\,dx\\
&=n\,i\int_{\R^n}\bar{u}\left(uv-\tfrac12 \Delta u \right)dx-\int_{\R^n}(x\cdot \nabla \bar{u})u_t\,dx\\
&=n\,i\int_{\R^n}|u|^2v\,dx -i\frac{n}{2}\int_{\R^n}\bar{u}\Delta u\,dx -\int_{\R^n}(x\cdot \nabla \bar{u})u_t\,dx\\
&=n\,i\int_{\R^n}|u|^2v\,dx +i\frac{n}{2}\int_{\R^n}|\nabla u|^2\,dx -\int_{\R^n}(x\cdot \nabla \bar{u})u_t\,dx.
\end{split}
\end{equation}
Combining (\ref{proof-virial-a}) and (\ref{proof-virial-b}) it follows that
\begin{equation}\label{proof-virial-c}
\frac{d}{dt}\int_{\R^n}(x\cdot \nabla u)\bar{u}\,dx=n\, i\int_{\R^n}v |u|^2dx +i\frac{n}{2}\int_{\R^n}|\nabla u|^2dx +
\int_{\R^n} x\cdot( \bar{u}_t\nabla u - u_t\nabla \bar{u})dx.
\end{equation}
On the other hand, using the first equation of the system and integrating by parts we get
\begin{equation}\label{proof-virial-d}
\begin{split}
\int_{\R^n} x\cdot(\bar{u}_t\nabla u - u_t\nabla \bar{u})dx
&=2 i\,\text{Im}\int_{\R^n}(x\cdot \nabla u)\bar{u}_t\, dx\\
&=2 i\,\text{Im}\int_{\R^n}i\,(x\cdot \nabla u)\left(\bar{u}v-\tfrac{1}{2}\Delta \bar{u}\right)dx\\
&=i\left(2\text{Re}\int_{\R^n}(x\cdot \nabla u)\bar{u}v\, dx - \text{Re}\int_{\R^n}(x\cdot \nabla u)\Delta \bar{u}\, dx \right)
\end{split}
\end{equation}
and
\begin{equation}\label{proof-virial-e}
\begin{split}
\text{Re}\int_{\R^n}(x\cdot \nabla u)&\Delta \bar{u}\, dx = \text{Re}\int_{\R^n}\sum\limits_j x_j\p_{x_j}u \sum\limits_k\p^2_{x_k}\bar{u}\,dx\\
&=\text{Re}\Biggl(-\int_{\R^n}|\nabla u|^2 dx - \int_{\R^n} \sum\limits_{j,k} x_j \p^2_{x_kx_j}u\p_{x_k}\bar{u}\,dx \Biggl)\\
&=-\int_{\R^n}|\nabla u|^2 dx -
\frac{1}{2}\Biggl(\int_{\R^n} \sum\limits_{j,k} x_j \p^2_{x_kx_j}u\p_{x_k}\bar{u}\,dx +
\int_{\R^n} \sum\limits_{j,k} x_j \p^2_{x_kx_j}\bar{u}\p_{x_k}u\,dx \Biggl)\\
&= -\int_{\R^n}|\nabla u|^2\,dx + \frac{n}{2}\int_{\R^n}|\nabla u|^2\,dx\\
&=\Bigl(\frac{n}{2}-1\Bigl)\int_{\R^n}|\nabla u|^2\,dx.
\end{split}
\end{equation}
Hence, from (\ref{proof-virial-d}) and (\ref{proof-virial-e}) we have

\begin{equation}\label{proof-virial-f}
\int_{\R^n} x\cdot(\bar{u}_t\nabla u - u_t\nabla \bar{u})dx = 2  i\, \text{Re}\int_{\R^n}(x\cdot \nabla u)\bar{u}v\,dx
-i\Bigl(\frac{n}{2}-1\Bigl)\int_{\R^n}|\nabla u|^2\,dx.
\end{equation}
Now, combining the last inequality with (\ref{proof-virial-c}) we obtain
\begin{equation}\label{proof-virial-g}
\begin{split}
\frac{d}{dt}\text{Im}\int_{\R^n}(x\cdot \nabla u)\bar{u}\,dx &= \int_{\R^n}|\nabla u|^2dx + n\int_{\R^n}v |u|^2dx +
2\, \text{Re}\int_{\R^n}(x\cdot \nabla u)\bar{u}v\, dx\\
&= \int_{\R^n}|\nabla u|^2dx + n\int_{\R^n}v |u|^2dx + \int_{\R^n}(x\cdot \nabla |u|^2) v\, dx.
\end{split}
\end{equation}
Finally, using (\ref{proof-virial-a0}) and the pseudo-Hamiltonian (\ref{Energy-1})  with $\lambda = -1$, we deduce the claimed
equality:
\begin{equation*}
\frac{d^2}{dt^2}\frac12 \int_{\R^n}|x|^2|u|^2dx = E(t)+(n-2)\int_{\R^n}v|u|^2\,dx -\int_{\R^n}v^2\,dx
+\int_{\R^n}( x\cdot \nabla |u|^2) v\,dx.
\end{equation*} 
As before, this procedure is 
made rigorous by following a regularization limiting method as in \cite{Cazenave-Book} or \cite{Merle}, so
this concludes the proof. $\square$

\section{\bf{Blow-up or Unboundedness of the Sobolev Norms of the Solutions for the Focusing Energy-Critical Case.}}\label{Section-blowup}

We begin this final section by providing a description of the method that will be pursued in proving Theorem \ref*{blowup-theorem}.
The proof follows by contradiction. Given the initial data in the statement of the theorem, the time $T_0$ is determined through an argument that is in itself a crucial part of the whole proof. 
We will leave it to the logical 
flow of the presentation to explain how that $T_0$ is obtained. But it will become clear that it is a function
only of the fixed initial data.
We thus begin by considering that such a $T_0$, whose dependence on the data will be established later, has been fixed and we assume the following two hypotheses:
\begin{enumerate}[(I)]
\item there is a $T \geq T_0$ and $0<\mu_0<1$ such that the solutions of the  Cauchy problem \eqref{SD}, 
given by Theorem \ref{LWP-4D}, for all values of the delay parameter
$0< \mu \le\mu_0$ and fixed initial data $(u_0,v_0) \in \mathscr{S}({\R}^4)\times \mathscr{S}({\R}^4)$,  exist on $[0,T]$, i.e. we assume that for any $0< \mu \le\mu_0$ the maximal time interval $(T_{\mu_{l_1}}, T_{\mu_{r_1}})$ of existence of the solution satisfies $T_0\leq T<T_{\mu_{r_1}}$; \smallskip 
\item 
there exists $s>4$ such that the $H^s$-norms of the solutions for the Schr\"odinger component $u_{\mu}$ of the system are uniformly bounded, in $\mu$, on the time interval $[0, T]$.
\end{enumerate}

\smallskip
By proving that these two hypotheses cannot both be true, through a contradiction argument, we will
therefore be proving that either (I) does not hold, or else (I) is true but (II) is not. And these
are precisely the two alternatives in the statement of Theorem \ref{blowup-theorem}. It should also be
noted that, although we will establish the existence of one such $T_0$, the optimal value should be considered 
for the final result. In other words, once we show that Theorem \ref{blowup-theorem} is true for one value
of $T_0>0$, the infimum of all such times, for which the theorem still holds, is the number to
be regarded as the $T_0$ in the statement of  Theorem \ref{blowup-theorem}.

\smallskip
Analogously to the convexity argument that is used for the classical blow-up results for the NLS equation, we will use the virial type identity of Theorem \ref{virial-theorem} to show that, under hypotheses (I) and (II), there exists a sequence of small delay parameters $\mu \to 0$ for which the variance \eqref{variance} of the corresponding solutions will all necessarily decrease to zero in $[0,T_0]$, thus leading to contradiction for the elements of this sequence. Recall also that, as pointed in  Remark \ref{persistence}, the solutions $(u,v)$ preserve the $H^s\times H^s$-regularity, for any $s>1$, in the  interval $(T_{\mu_{l_1}}, T_{\mu_{r_1}})$, so that $(T_{\mu_{l_1}}, T_{\mu_{r_1}})=(T_{\mu_{l_s}}, T_{\mu_{r_s}})$ and we can use \eqref{virial-b} because our initial data, in the Schwartz class, is infinitely regular. 

\smallskip

The general framework of the proof consists of splitting the analysis into two steps, corresponding to two fixed and consecutive time intervals
$[0,t_0]$ and $[t_0,T_0\,]$, independent of $\mu$. Unlike in the usual proof for the NLS equation, where the virial identity
can be used to employ the convexity argument for the time evolution of the variance \eqref{variance}, starting right at the initial time $t=0$, this cannot be done in our case.
In fact, formula \eqref{virial-b} is significantly more complicated than the analogous formula for the NLS equation
and its negativity as a consequence of the negative initial pseudo-energy $E_0<0$ can
only be attained through a careful
control of the nonlinear terms that result from the Debye relaxation, for a fixed interval $[t_0,T_0\,]$, with large enough $t_0$ and an adequate sequence of small parameters $\mu\to 0$. This is not surprising as only for large times and small parameter $\mu$ can we expect the system \eqref{SD} to behave similarly 
to the NLS equation, due to the Debye delay effect. 

\smallskip 
On the other hand, the first time interval $[0,t_0]$ satisfies a dual role; it postpones the starting time
of the convexity arguments to a big enough value of $t_0$, to be determined a priori, 
while also enabling the estimate of uniform bounds for the variance and its first derivative at $t_0$,
independently of (small) $\mu$. These will be used subsequently as the coefficients of a convex parabola, starting at $t_0$, which serves as an upper
bound for the time evolution for $t\geq t_0$ of the variance of all solutions, with small $\mu$, 
thus allowing the 
maximum required length of the second interval $[t_0,T_0\,]$ to be determined, in order
to achieve the contradiction for all such solutions. This is, therefore, the method for
obtaining $T_0$.

\smallskip 
At each of the two time intervals a limiting procedure, as $\mu \to 0$, will be employed to
yield a reference function
with which uniform estimates can be obtained independently of the small values of $\mu$; on $[0,t_0]$ we estimate upper bounds for the
variance and its first derivative at $t_0$,
while on $[t_0,T_0\,]$ we estimate the sign of \eqref{virial-b}. These estimates, as well as the limiting
procedure, can only be proved for fixed time intervals. Therefore, each time interval has to
be determined in advance of the corresponding step of the proof; $t_0$ is chosen at the beginning so large that the negative sign of $E_0$ will later become dominant in \eqref{virial-b} for $t \geq t_0$, while
$T_0$ is estimated afterwards as a foreseeable upper bound for the contradiction to occur, uniformly for small
$\mu$, due to the convex time evolution of the variances starting at $t_0$.  

\smallskip 
To finish this introductory description of the strategy that we will pursue, we point out that, actually, the computations and
mathematical techniques used in both steps of the proof are very similar. In either case, on $[0,t_0]$ and on $[t_0,T_0\,]$,
the plan is ultimately to estimate the second derivative of the variance \eqref{virial-b} which, after two integrations, provides an upper bound for the time evolution of the variance itself. The only difference
resides in the fact that, on $[0,t_0]$ we intend to obtain estimates for the variance and its first derivative at the
end of the interval $t_0$, whereas on $[t_0, T_0\,]$ the goal is to show that the time evolution of the variance
is convex, i.e. with negative second derivative, so that it becomes zero within that interval.

\subsection{\bf{Useful results.}}For the purpose of applying the limiting procedure in a straightforward manner during the proof, we now establish the
following results.

\smallskip

\begin{lemma}\label{lemauniforme}
	Let $(u_0,v_0)$ be radially symmetric functions belonging to the space $H^{s+2}(\R^n) \times H^{s+2}(\R^n)$, with $s>n/2$,  and suppose that for all  $0< \mu \le \mu_0$, for some $\mu_0 >0$, the
	solutions $u_{\mu}$ of the integro-differential equation \eqref{SD-IDF}, given by Theorem 2 in \cite{Bidegaray2},
	are defined on a common fixed time interval $[0,T]$, i.e. $u_{\mu} \in C([0,T];H^{s+2}(\R^n))$, $\forall_{0<\mu\leq\mu_0}$,  satisfying 
	$$\sup\limits_{0 < \mu \le \mu_0} \Big(\sup\limits_{0\le t \le T}\|u_{\mu}(\cdot, t)\|_{H^s}\Big)=K_s.$$
	Then, there exist  positive constants $M_s$ and $N_s$, independent of $\mu$, such that for all $0< \mu \le \mu_0$, the following uniform bounds also hold 
	\begin{enumerate}[(a)]
		\item $\displaystyle \sup_{[0,T]} \|\partial_t u_{\mu}(\cdot, t)\|_{H^s} \leq M_s,$ \smallskip
		\item $\displaystyle \sup_{[0,T]} \| |x| u_{\mu}(\cdot, t)\|_{L^{\infty}}\le N_s.$
	\end{enumerate}
\end{lemma}

\begin{remark}
It should be noted that the $u$ component of the solution  obtained in Theorem \ref{LWP-4D}, by performing a fixed point argument to the full
system \eqref{SD}, is always a solution of the integro-differential equation \eqref{SD-IDF}. For high Sobolev regularity $s>n/2$ the converse is also
easily seen to be true, i.e. that solutions $u$ of the integro-differential equation yield solutions of the system by simply defining $v$ to be given by \eqref{equation-v}.
Furthermore, we note that the radial hypothesis on the data is not necessary for the conclusion in (a),  as can be seen in the proof below.
\end{remark}

\noindent \textbf{Proof of Lemma \ref{lemauniforme}.} Using the algebra properties of $H^s$ for $s>n/2$ we estimate $u_t$ directly from \eqref{SD-IDF} in the following way
$$
\|\partial_tu_{\mu}\|_{H^{s}} \leq \frac12\|u_{\mu}(t)\|_{H^{s+2}} + \|v_0\|_{H^s}\|u_{\mu}(t)\|_{H^s}+
\|u_{\mu}(t)\|_{H^s}\frac{1}{\mu}\int_0^{t}e^{-\frac{t-\tau}{\mu}}\|u_{\mu}(\tau)\|^2_{H^s} d\tau,
$$
from which we get
$$
\|\partial_tu_{\mu}\|_{L^{\infty}_{[0,T]}H^{s}}\leq \frac12 K_s + \|v_0\|_{H^{s+2}}K_s
+ T K_s^3,
$$
again depending only on the initial data, $s$ and the time interval $[0,T]$, but independent of $\mu$. And that
is (a).

\smallskip
On the other hand,  for radial initial data the solutions remain radial for all times, so that for high enough regularity the radial version of the Gagliardo-Nirenberg inequality, together with Sobolev embedding, enable the control
of the norm $\||x| u_{\mu}\|_{L^{\infty}}$ (see Radial Lemma in pp. 155 of  \cite{Strauss}). 
$\square$

\smallskip

\begin{corollary}\label{limitfunction}
	Under the hypotheses of the previous lemma, for data $(u_0,v_0) \in H^{s+2}(\R^n) \times H^{s+2}(\R^n)$,
	with $s>n/2$,
	there exists a function $\tilde{u} \in C([0,T]; H^{s}(\R^n))$ and a sequence $u_{\mu_i}$ with
	$\mu_i \to 0$ such that
	$$\lim_{i} \|u_{\mu_i} - \tilde{u}\|_{L^{\infty}_{[0,T]}H^s_x}=0.$$
\end{corollary}

\noindent \textbf{Proof.}  From the previous lemma, this is a direct application of the Arzel\`a-Ascoli Theorem 
for the family of functions $u_{\mu}$ over the time interval $[0,T]$: the hypothesis of the lemma  provides uniform boundedness in $H^{s+2}$, and therefore also in $H^s$, while condition (a) provides equicontinuity in $H^s$.
$\square$

\subsection{\bf{Proof of Theorem \ref*{blowup-theorem}.}}

We start by defining the nonnegative function $h$, as (half of) the variance of $u$,
$$h(t)=\displaystyle \frac12 \int_{\R^n} |x|^2|u|^2dx.$$
As we are assuming that, $0<\mu\leq \mu_0$, the corresponding solution $u=u_{\mu}$ exists on $[0,T]$, 
then $h$ must be strictly positive for all $t\in[0,T]$, because the hypothesis $E_0<0$ 
guarantees that $u_0 \not= 0$ so that $h(0)>0$, and if there existed a time for which $h$ would 
vanish, then, from Heisenberg's inequality and the $L^2$ conservation
of $u$, we would get
$$0<\|u_0\|_{L^2}^2=\|u(\cdot,t)\|^2_{L^2}\leq \frac{1}{2} \| |x| u(\cdot,t)\|_{L^2}\|\nabla u(\cdot,t)\|_{L^2},$$
so that $\|\nabla u(\cdot,t)\|_{L^2}$ would blow up at that same point, contradicting the existence of the solution up to time $T$.

\smallskip 

We now proceed to determine the length of the first time interval, i.e. $t_0$, which will also be
the starting time for the convexity arguments based on the virial identity, uniformly for small $\mu \to 0$. As $E_0 <0$, and using the $L^2$ conservation property of the solution $u$ \eqref{Conservation Law}, 
we can conclude that there exists a large enough positive time, $t_0>0$, for which
\begin{equation}\label{novidade}
\begin{split}
\left| e^{-t/\mu}\!\!\int_{\R^4}\!\!(2v_0 + x\cdot \nabla v_0)|u|^2\,dx\right| 
&\leq e^{-t_0} \|2v_0 + x\cdot \nabla v_0\|_{L^{\infty}_x} \|u\|_{L^2_x}^2 \\
&= e^{-t_0} \|2v_0 + x\cdot \nabla v_0\|_{L^{\infty}_x} \|u_0\|_{L^2_x}^2\\
&\leq \frac{|E_0|}{2},
\end{split}
\end{equation}
for all $t \geq t_0$ and $0<\mu \leq\mu_0$. We will see later, in the second part of the proof, that this condition will guarantee
the existence of a negative upper bound for $h''$ over the second interval $[t_0,T_0\,]$.

\smallskip
With the first time interval $[0,t_0]$ determined and fixed, we can now move on to estimating $h''$
in order to obtain, after two integrations over the length of the interval, uniform estimates for $h(t_0)$ and $h'(t_0)$, for small values of $\mu \to 0$.

\smallskip
Using, from \eqref{Energy-1}
and $\lambda=-1$, the fact that
$E(t) \le E_0$ in (\ref{virial-b}), we obtain
\begin{equation}\label{proof-blowup-1}
\begin{split}
h''(t)&\le E_0 +(n-2)\int_{\R^n}v|u|^2\,dx-\int_{\R^n}v^2\,dx +\int_{\R^n} (x \cdot \nabla |u|^2) v\,dx
\\
& \leq  E_0 +(n-2)\int_{\R^n}v|u|^2\,dx+\int_{\R^n} (x \cdot \nabla |u|^2) v\,dx, 
\end{split}
\end{equation}
for all $t\geq 0$.

\smallskip
Applying now the integral form (\ref{equation-v}) of the $v$-solution in the previous inequality we obtain
\begin{equation}\label{proof-blowup-2}
\begin{split}
h''(t)&
\le E_0+(n-2)\left(e^{-t/\mu}\int_{\R^n}v_0|u|^2\,dx-\frac1{\mu}\int_0^te^{-\frac{t-\tau}{\mu}}\int_{\R^n}|u(t)|^2|u(\tau)|^2\,dxd\tau\right)\\
&\hspace{1.2cm}+e^{-t/\mu}\int_{\R^n}v_0( x \cdot \nabla |u|^2)\,dx
-\frac{1}{\mu}\int_0^t e^{-\frac{t-\tau}{\mu}}\int_{\R^n}(x\cdot \nabla |u(t)|^2)|u(\tau)|^2\,dxd\tau
\end{split}
\end{equation}
and we rewrite the right hand side of (\ref{proof-blowup-2}) as follows:
\begin{equation}\label{proof-blowup-2-a}
\begin{split}
h''(t)&\le E_0+e^{-t/\mu}\int_{\R^n}\left( (n-2)v_0 + v_0( x \cdot \nabla |u|^2)\right) \,dx
-\frac{n-2}{\mu}\int_0^te^{-\frac{t-\tau}{\mu}}\int_{\R^n}|u(t)|^4\,dx d\tau\\
&\hspace{1.2cm}-\frac{1}{\mu}\int_0^te^{-\frac{t-\tau}{\mu}}\int_{\R^n}(x\cdot \nabla |u(t)|^2)|u(t)|^2\,dxd\tau
+r(\mu, t),
\end{split}
\end{equation}
where $r(\mu, t)=r_1(\mu, t)+ r_2(\mu, t)$ with
\begin{align}
&r_1(\mu, t)=\frac{2-n}{\mu}
\int_0^te^{-\frac{t-\tau}{\mu}}\int_{\R^n}|u(t)|^2\left(|u(\tau)|^2- |u(t)|^2\right)
\,dxd\tau,\label{proof-blowup-2-b}\\
&r_2(\mu, t)=-\frac{1}{\mu}
\int_0^te^{-\frac{t-\tau}{\mu}}\int_{\R^n}(x\cdot \nabla |u(t)|^2)\left(|u(\tau)|^2- |u(t)|^2\right)\,dxd\tau\label{proof-blowup-2-c}.
\end{align}

\smallskip 
Integrating by parts the second and fourth terms in (\ref{proof-blowup-2-a}), we get
\begin{equation}\label{proof-blowup-3}
\begin{split}
h''(t)&\le E_0-e^{-t/\mu}\int_{\R^n}(2v_0 + x\cdot \nabla v_0)|u|^2\,dx
-\frac{n-2}{\mu}\int_0^te^{-\frac{t-\tau}{\mu}}\int_{\R^n}|u(t)|^4\,dxd\tau\\
&\hspace{1.2cm} +\frac{n}{2\mu}\int_0^te^{-\frac{t-\tau}{\mu}}\int_{\R^n}|u(t)|^4\,dxd\tau
+r(\mu, t)\\
&=E_0-e^{-t/\mu}\!\!\int_{\R^n}\!\!(2v_0 + x\cdot \nabla v_0)|u|^2\,dx +
\left(2-\frac{n}{2}\right)\left(1-e^{-t/\mu}\right)\|u(t)\|_{L^4_x}^4\\
&\hspace{1.2cm}+r(\mu, t).
\end{split}
\end{equation}
Up until this point, the computations were performed for any number of spatial dimensions and all $t \geq 0$. 
From here on we 
restrict to the $n=4$ case
for which the theorem is stated. The third term on the right hand side of the previous inequality disappears
\begin{equation}\label{upperboundh}
h''(t) \leq E_0-e^{-t/\mu}\!\!\int_{\R^n}\!\!(2v_0 + x\cdot \nabla v_0)|u|^2\,dx +
r(\mu,t),
\end{equation}
and we are left with estimating $r(\mu, t)$.

\smallskip
We rewrite (\ref{proof-blowup-2-b})-(\ref{proof-blowup-2-c}) as follows:
\begin{equation}\label{proof-blowup-4}
r_j(\mu, t)=\frac{1}{\mu}\int_0^te^{-\frac{(t-\tau)}{\mu}}I_j(t,\tau)d\tau
\quad (j=1,2),
\end{equation}
where
\begin{align}
&I_1(t, \tau)=-2\int_{\R^4}|u(t)|^2\left(|u(\tau)|^2- |u(t)|^2\right)\,dx,\label{proof-blowup-5-a}\\
&I_2(t, \tau)=-\int_{\R^4}(x\cdot \nabla |u(t)|^2)\left(|u(\tau)|^2- |u(t)|^2\right)\,dx\label{proof-blowup-5-b}.
\end{align}

\smallskip
Next, we estimate the integrals $I_1(t, \tau)$ and $I_2(t, \tau)$. More precisely, using H\"older's inequality and (\ref{Gagliardo-Nirenberg}) we have
\begin{equation}\label{proof-blowup-6-a}
\begin{split}
|I_1(t, \tau)|&\le 2\|u(t)\|^2_{L^4_x}\|\,|u(\tau)|+|u(t)|\,\|_{L^4_x}\|\,|u(\tau)|-|u(t)|\, \|_{L^4_x}\\
&\le c \|\nabla u(t)\|^2_{L^2_x}\left(\|\nabla u(\tau)\|_{L^2_x}+ \|\nabla u(t)\|_{L^2_x}\right)\|\nabla u(\tau) -\nabla u(t)\|_{L^2_x}\\
&\le c  \|u\|^3_{L^{\infty}_{[0,t_0]}H^1_x}\|\nabla u(\tau) -\nabla u(t)\|_{L^2_x},
\end{split}
\end{equation}
for all $0\le \tau \le t\le t_0$, so that
\begin{equation}\label{proof-blowup-8-a}
\begin{split}
|r_1(t,\mu)|&\le \frac{1}{\mu}\int_0^te^{-\frac{(t-\tau)}{\mu}}|I_1(t,\tau)|d\tau\\
&\le c\|u\|^3_{L^{\infty}_{[0,t_0]}H^1_x}\frac{1}{\mu}\int_0^te^{-\frac{(t-\tau)}{\mu}}\|\nabla u(\tau) -\nabla u(t)\|_{L^2_x}d\tau.
\end{split}
\end{equation}
For the first time interval $[0,t_0]$ we are not concerned with obtaining a very fine estimate for $h''$, as
we only wish to get uniform bounds for $h(t_0)$ and $h'(t_0)$. Therefore, we simply do
$$\|\nabla u(\tau) -\nabla u(t)\|_{L^2_x} \leq 2 \|u\|_{L^{\infty}_{[0,t_0]}H^1_x},$$
and thus
\begin{equation}\label{R1int1}
\begin{split}
|r_1(t,\mu)|&\le c\|u\|^4_{L^{\infty}_{[0,t_0]}H^1_x}\frac{1}{\mu}\int_0^te^{-\frac{(t-\tau)}{\mu}}d\tau\\
&\le c\|u\|^4_{L^{\infty}_{[0,t_0]}H^1_x}.
\end{split}
\end{equation}

\smallskip
In a similar way we have
\begin{equation}\label{proof-blowup-6-b}
|I_2(t, \tau)|\le c \|u\|^2_{L^{\infty}_{[0,t_0]}H^1_x} \||x| u\|_{L^{\infty}_{[0,t_0]}L^{\infty}_x}
\|\nabla u(\tau) -\nabla u(t)\|_{L^2_x},
\end{equation}
for all $0\le \tau \le t\le t_0$ and
\begin{equation}\label{proof-blowup-8-b}
\begin{split}
|r_2(t,\mu)|&\le \frac{1}{\mu}\int_0^te^{-\frac{(t-\tau)}{\mu}}|I_1(t,\tau)|d\tau\\
&\le c\|u\|^2_{L^{\infty}_{[0,t_0]}H^1_x} \||x| u\|_{L^{\infty}_{[0,t_0]}L^{\infty}_x}
\frac{1}{\mu}\int_0^te^{-\frac{(t-\tau)}{\mu}}\|\nabla u(\tau) -\nabla u(t)\|_{L^2_x}d\tau.
\end{split}
\end{equation}
Again, doing the simple estimate as was performed before for $r_1(t,\mu)$, we get
\begin{equation}\label{R2int1}
|r_2(t,\mu)|
\le c\|u\|^3_{L^{\infty}_{[0,t_0]}H^1_x} \||x| u\|_{L^{\infty}_{[0,t_0]}L^{\infty}_x}.
\end{equation}

\smallskip
At this point we appeal to the limiting procedure of Corollary \ref{limitfunction}, to 
 yield a reference function $\tilde{u}$ with respect
 to which uniform estimates, independent of $\mu$, can be obtained. As our initial data is in the Schwartz space, an
 arbitrarily high Sobolev regularity can be chosen in order to apply Corollary \ref{limitfunction} to
 conclude that there is a converging sequence of solutions $u=u_{\mu_i}$, with $\mu_i \to 0$, in the
 $L^{\infty}_{[0,t_0]}H^1_x$ norm. Their norms are therefore bounded, while the $L^{\infty}_{[0,t_0]}L^{\infty}_x$
 norm is uniformly bounded in $\mu$, from Lemma \ref{lemauniforme}.
 
 \smallskip
 We can thus conclude, from \eqref{R1int1} and \eqref{R2int1}, that the uniform boundedness of the norms
 implies that there exists a constant $R$ such that, for all solutions in the sequence above $u=u_{\mu_i}$,
 we have
 $$r(t,\mu_i)=r_1(t,\mu_i)+r_2(t,\mu_i) \leq R,$$
 for all $t \in [0,t_0]$ and independently of $\mu_i$.
 
 \smallskip
 This uniform bound for $r(t,\mu_i)$ and the estimate \eqref{upperboundh} now imply that, for $t \in [0,t_0]$, we have
 $$h''(t) \leq E_0-e^{-t/\mu}\!\!\int_{\R^4}\!\!(2v_0 + x\cdot \nabla v_0)|u|^2\,dx + R,$$
 from which a uniform bound can be easily obtained
 \begin{equation}
 \begin{split}
 h''(t)&\leq E_0+\left| e^{-t/\mu}\!\!\int_{\R^4}\!\!(2v_0 + x\cdot \nabla v_0)|u|^2\,dx\right|  + R\\
 &\leq E_0 +  \|2v_0 + x\cdot \nabla v_0\|_{L^{\infty}_x} \|u_0\|_{L^2_x}^2 + R,
 \end{split}
 \end{equation}
 that only depends on the initial data and the constant $R$, but not on $\mu_i$.
 
 \smallskip
 To finish the first half of the proof of the theorem, we just need to integrate the previous inequality twice, from
 $0$ to $t_0$, 
 $$h'(t_0)\leq h'(0) + \left(E_0 +  \|2v_0 + x\cdot \nabla v_0\|_{L^{\infty}_x} \|u_0\|_{L^2_x}^2 + R\right) t_0,$$
 and
 $$h(t_0) \leq h(0) + h'(0) t_0 + \left(E_0 +  \|2v_0 + x\cdot \nabla v_0\|_{L^{\infty}_x} \|u_0\|_{L^2_x}^2 + R\right)
 \frac{t_0^2}{2},$$
 recalling also that
 $$h(0)=\frac12 \int_{\R^n} |x|^2|u_0|^2dx \qquad \mbox{and} \qquad h'(0)=\emph{Im}\int_{\R^n}(x\cdot \nabla u_0)\overline{u_0}\,dx,$$
 to
 conclude that there exist two constants, let us call them $A,B >0$, depending only
 on the initial data and $t_0$, such that
 \begin{equation}\label{ABbounds}
 h(t_0)\leq A \qquad \mbox{and}\qquad h'(t_0)\leq B,
 \end{equation}
 uniformly, for all the solutions in the previous sequence $u=u_{\mu_i}$ with $\mu \to 0$.
 
 \bigskip
 We now start the second half of the proof, by first determining the length of the second time
 interval $[t_0,T_0\,]$, and consequently the crucial value of $T_0$. This interval is where the convexity of the time
 evolution of the variance $h$, for all solutions $u=u_{\mu_i}$, starting at $t_0$ will lead to the contradiction. We
 will eventually prove that, on this interval, $h''(t) \leq E_0 /4<0$ uniformly in $\mu_i$, which,
 together with the bounds \eqref{ABbounds} implies that the convex parabola
 \begin{equation}\label{parabola}
 g(t)=\frac{E_0}{8}(t-t_0)^2 + B (t-t_0)+A,
 \end{equation}
 is a uniform upper bound for the time evolution of all variances $h(t)$ for $t \geq t_0$.
Therefore, the first root of this parabola larger than $t_0$ (which is well defined, as $A>0$) is an upper bound
in time for all variances $h(t)$ to become zero, for all solutions $u=u_{\mu_i}$. We thus choose $T_0$
as the first root of $g$, which we can see now is only dependent on the initial data and could have been 
computed right at the beginning of the proof, as required.
 
 \smallskip
 With the second time interval $[t_0,T_0\,]$ determined and fixed, we proceed to estimate
 $h''$ on it. Compared to the first half of the proof, on $[0,t_0]$, the goal now is
 to actually establish a uniform negative upper bound for this second derivative in time, which demands slightly finer estimates than before.
  The computations from \eqref{proof-blowup-1} to \eqref{proof-blowup-8-b} are repeated in exactly
 the same way as in the first part, except that now the integrations in time are performed from
 $0$ to $t \in [t_0,T_0\,]$, while the $L^{\infty}$ norms in time are taken over the whole
 time interval $[0,T_0\,]$.
 As before, the final step consists in estimating the reminder term $r(t,\mu) =r_1(t,\mu)+ r_1(t,\mu)$.
 But unlike in the first interval, where we only needed to obtain uniform bounds for these terms,
 now, in the second interval, we actually need to show that they can be made arbitrarily small
 uniformly in $\mu_i$, so that the negativity of $E_0$ dominates the upper bound of $h''$ in 
 \eqref{upperboundh}, over the whole interval $[t_0,T_0\,]$.
 
 \smallskip
 To clarify the presentation, we gather here again the two formulas for the reminder terms that
will be handled at this point,
\begin{equation}\label{remind1}
|r_1(t,\mu_i)|\le c\|u\|^3_{L^{\infty}_{[0,T_0\,]}H^1_x}\frac{1}{\mu_i}\int_0^te^{-\frac{(t-\tau)}{\mu_i}}\|\nabla u(\tau) -\nabla u(t)\|_{L^2_x}d\tau,
\end{equation}
and
\begin{equation}\label{remind2}
|r_2(t,\mu_i)|\le c\|u\|^2_{L^{\infty}_{[0,T_0\,]}H^1_x} \||x| u\|_{L^{\infty}_{[0,T_0\,]}L^{\infty}_x}
\frac{1}{\mu_i}\int_0^te^{-\frac{(t-\tau)}{\mu_i}}\|\nabla u(\tau) -\nabla u(t)\|_{L^2_x}d\tau,
\end{equation}
for $t_0 \leq t \leq T_0.$

\smallskip
We will now show that the integral
\begin{equation}\label{proof-blowup-8-c}
\frac{1}{\mu_i}\int_0^te^{-\frac{(t-\tau)}{\mu_i}}\|\nabla u(\tau) -\nabla u(t)\|_{L^2_x}d\tau,
\end{equation}
present in both \eqref{remind1} and \eqref{remind2}, can be made arbitrarily small, uniformly for 
all $t_0\leq t \leq T_0$, if
we pick any $\mu_i$ smaller than a conveniently chosen $\tilde{\mu}_0$, that depends only on $t_0$, $T_0$ and on the desired
smallness of \eqref{proof-blowup-8-c}. For that we will exploit the fact that the function $\frac{1}{\mu}e^{-\frac{t}{\mu}}$, for $t>0$, is essentially an
approximate identity as $\mu \to 0$. However, $u$ also depends on $\mu$ which prevents a direct approach to show
that $r(t,\mu_i)$ is small as $\mu_i \to 0$. So, again at this point, we use the limiting procedure of
Corollary \ref{limitfunction} in order to yield a reference function $\tilde{u}$ on which the
approximate identity can be applied and estimates can be uniformly obtained. 

\smallskip
So, applying Corollary \ref{limitfunction} to the sequence of solutions $u=u_{\mu_i}$ inherited from the first
part of the proof, and now considering the full time interval $[0,T_0\,]$, we can extract
a subsequence - which we will continue denoting by $u_{\mu_i}$ - and a limit function
$\tilde{u} \in C([0,T_0\,];H^1(\R^4))$ such that $\displaystyle \lim_{\mu_i \to 0} \|u_{\mu_i} - \tilde{u}\|_{L^{\infty}_{[0,T_0\;]}H^1_x}=0.$

\smallskip
Given now any arbitrary $\varepsilon >0$, using the continuity of $\tilde{u}$, as a flow from the closed time interval $[0, T_0\,]$ to $H^1$, we can take a positive number $\delta_{\varepsilon, T_0}$ such that
\begin{equation}\label{proof-blowup-7}
\|\nabla \tilde{u}(\tau) -\nabla \tilde{u}(t)\|_{L^2_x}\le \varepsilon,\;\; \text{for all}\;\;0\leq \tau \leq t \leq T_0\;\;
\text{with}\;\;t -\tau \le \delta_{\varepsilon, T_0}.
\end{equation}

\smallskip
Without loss of generality, we can obviously assume that $\delta_{\varepsilon, T_0}<t_0$ so
that $t-\delta_{\varepsilon, T_0}>0$, for $t_0 \leq t \leq T_0$, and we can break the integral \eqref{proof-blowup-8-c} in two,
\begin{equation}\label{proof-blowup-9}
\frac{1}{\mu_i}\int_0^{t-\delta_{\varepsilon, T_0}} e^{-\frac{(t-\tau)}{\mu_i}}\|\nabla \tilde{u}(\tau) -\nabla \tilde{u}(t)\|_{L^2_x}d\tau
+
\frac{1}{\mu_i}\int_{t-\delta_{\varepsilon, T_0}}^t e^{-\frac{(t-\tau)}{\mu_i}}\|\nabla \tilde{u}(\tau) -\nabla \tilde{u}(t)\|_{L^2_x}d\tau.
\end{equation}
We estimate the first of these two
\begin{equation}\label{proof-blowup-9-a}
\begin{split}
\frac{1}{\mu_i}\int_0^{t-\delta_{\varepsilon, T_0}} 
e^{-\frac{(t-\tau)}{\mu_i}}\|\nabla \tilde{u}(\tau) -\nabla \tilde{u}(t)\|_{L^2_x}d\tau 
&\leq 2\|\nabla \tilde{u}\|_{L^{\infty}_{[0,T_0\,]}L^2_x}\,\,\frac{1}{\mu_i}\int_0^{t-\delta_{\varepsilon, T_0}} 
e^{-\frac{(t-\tau)}{\mu_i}}d\tau \\
&= 2\|\nabla \tilde{u}\|_{L^{\infty}_{[0,T_0\,]}L^2_x} \left( e^{-\frac{\delta_{\varepsilon, T_0}}{\mu_i}}-
e^{-\frac{t}{\mu_i}}\right)\\
&\leq \varepsilon\,\|\nabla \tilde{u}\|_{L^{\infty}_{[0,T_0\,]}L^2_x},
\end{split}
\end{equation}
uniformly for all $t_0\leq t\leq T_0$ and $0< \mu_i \leq\tilde{\mu}_0$, by
choosing a conveniently small $\tilde{\mu}_0=\tilde{\mu}_0(\varepsilon, \delta_{\varepsilon, T_0\,})$. As for the second integral
in \eqref{proof-blowup-9}, we use \eqref{proof-blowup-7} to estimate
\begin{equation}\label{proof-blowup-9-b}
\frac{1}{\mu_i}\int_{t-\delta_{\varepsilon, T_0}}^t
e^{-\frac{(t-\tau)}{\mu_i}}\|\nabla \tilde{u}(\tau) -\nabla \tilde{u}(t)\|_{L^2_x}d\tau 
\leq \varepsilon\,\,\frac{1}{\mu_i}\int_{t-\delta_{\varepsilon, T_0}}^t
e^{-\frac{(t-\tau)}{\mu_i}}d\tau \leq \varepsilon,
\end{equation}
also uniformly in $t \in [t_0,T_0\,]$ and in this case independently of $\mu_i$. Gathering both integrals,
we conclude that, for all $t_0\leq t \leq T_0$ and any $0<\mu_i\leq \tilde{\mu}_0$, we have
\begin{equation}\label{proof-blowup-9-c}
\frac{1}{\mu_i}\int_0^t
e^{-\frac{(t-\tau)}{\mu_i}}\|\nabla \tilde{u}(\tau) -\nabla \tilde{u}(t)\|_{L^2_x}d\tau 
\leq \varepsilon(1+\|\nabla \tilde{u}\|_{L^{\infty}_{[0,T_0\,]}L^2_x}).
\end{equation}

\smallskip
At this point, this estimate needs to passed on uniformly to the elements of the sequence $u=u_{\mu_i}$, which are
the ones that actually appear in \eqref{proof-blowup-8-c}. But this is easily done as a consequence of the
limit $\lim_i u_{\mu_i}=\tilde{u}$ in the $L^{\infty}_{[0,T_0\;]}H^1_x$ norm. We thus have
\begin{equation}\label{aproxima}
\begin{split}
\frac{1}{\mu_i}\int_0^te^{-\frac{(t-\tau)}{\mu_i}}\|\nabla u_{\mu_i}(\tau) -\nabla u_{\mu_i}(t)\|_{L^2_x}d\tau
&\leq 2\|\nabla u_{\mu_i} -\nabla \tilde{u}\|_{L^{\infty}_{[0,T_0\,]}L^2_x}\,\,\frac{1}{\mu_i}\int_0^{t} 
e^{-\frac{(t-\tau)}{\mu_i}}d\tau \\
& \hspace{2cm}+\frac{1}{\mu_i}\int_0^t
e^{-\frac{(t-\tau)}{\mu_i}}\|\nabla \tilde{u}(\tau) -\nabla \tilde{u}(t)\|_{L^2_x}d\tau \\
&\leq 2\|u_{\mu_i} - \tilde{u}\|_{L^{\infty}_{[0,T_0\;]}H^1_x} + \varepsilon(1+\|\nabla \tilde{u}\|_{L^{\infty}_{[0,T_0\,]}L^2_x}).
\end{split}
\end{equation}
We conclude, finally, that by making possibly $\tilde{\mu}_0$ even smaller, the term $\|u_{\mu_i} - \tilde{u}\|_{L^{\infty}_{[0,T_0\;]}H^1_x}$ can be also be made as small as desired, uniformly in $\mu_i \leq \tilde{\mu}_0$.
\smallskip
Hence, for any given small $\epsilon >0$, the estimates \eqref{remind1} and \eqref{remind2} become
\begin{equation}\label{proof-blowup-10}
|r_1(t,\mu_i)|\le c\,\epsilon\,\|u\|^3_{L^{\infty}_{[0,T_0\,]}H^1_x},
\end{equation}
and
\begin{equation}\label{proof-blowup-11}
|r_2(t,\mu_i)|\le c\,\epsilon\,\|u\|^2_{L^{\infty}_{[0,T_0\,]}H^1_x} \||x| u\|_{L^{\infty}_{[0,T_0\,]}L^{\infty}_x},
\end{equation}
for all $t\in [t_0,T_0\,]$ and  $0< \mu_i \le \tilde{\mu}_0(\varepsilon, T_0)$. But
the sequence $u=u_{\mu_i}$ is convergent in the $L^{\infty}_{[0,T_0\,]}H^1_x$ norm, so the
corresponding norms above are bounded, while the $\||x| u\|_{L^{\infty}_{[0,T_0\,]}L^{\infty}_x}$ is
also uniformly bounded in $\mu$, from Lemma \ref{lemauniforme}.

\smallskip 
So, for arbitrarily small $\eta$, choosing a suitable $\epsilon$ in (\ref{proof-blowup-10}) and (\ref{proof-blowup-11}), we can then take $\tilde{\mu}_0$ small enough so that
$$|r(t,\mu)|\le |r_1(t,\mu)|+|r_2(t,\mu)|\le \eta,$$
for all $t\in [t_0, T_0\,]$ and  $0< \mu \le \tilde{\mu}_0$. Thus, from \eqref{upperboundh}, for 
$t\in [t_0,T_0\,]$, we have
\begin{equation}\label{proof-blowup-12}
h''(t) \le E_0-e^{-t/\mu}\!\!\int_{\R^4}\!\!(2v_0 + x\cdot \nabla v_0)|u|^2\,dx + \eta.
\end{equation}

\smallskip
Condition \eqref{novidade}, which was pivotal in the choice of the length of the first interval, $t_0$,
finally makes its appearance, revealing the reason for that seemingly awkward option at the beginning of the proof: for $t \geq t_0$ the first two terms on the right hand side of \eqref{proof-blowup-12} are smaller
than $E_0/2$ so that
$$
h''(t) \leq \frac{E_0}{2} + \eta,
$$
and we can choose $\eta=|E_0|/4$ to get
\begin{equation}
h''(t) \leq \frac{E_0}{4},
\end{equation} 
for all $t \in [t_0,T_0\,]$ and $\mu_i\le \tilde{\mu}_0$ (recall that $\tilde{\mu}_0\leq \mu_0 <1$). 

\smallskip 
Consequently, $h(t)\le g(t)$ for $t\in [t_0,T_0\,]$ and we conclude  that the functions $h(t)$,
for every $\mu_i$ in the sequence, would necessarily  become zero, at some time instant in $[t_0,T_0\,]$, which is in contradiction 
with assumption (I), that for all $0<\mu\leq\mu_0$, the solutions $u_{\mu}$ exist for all $t \in [0,T]$.

\smallskip
Hypotheses (I) and (II) cannot therefore be simultaneously true, and this finishes the proof of the existence
of the two alternatives in Theorem \ref{blowup-theorem}. 

\smallskip
In case blow-up of the ${H^1(\R^4)\times  H^1(\R^4)}$ solutions does occur in finite time $t^*$, be it for
$t^* \leq T_0$ in alternative (a), or possibly also in alternative (b) for $t^* > T_0$, then
this implies that
$$\lim\limits_{t\nearrow t^*}\|(u,v)\|_{H^1(\R^4)\times  H^1(\R^4)}=+\infty.$$
From the $L^2$ conservation law of $u$ \eqref{Conservation Law} it follows that we must then have
\begin{equation}\label{blowupfinal}
\lim\limits_{t\nearrow t^*} (\|\nabla u\|_{L^2_x}+\|\nabla v\|_{L^2_x}+\|v\|_{L^2_x})=+\infty.
\end{equation}

% In fact, there would be a time $t^*\leq t_0$ for which $h(t^*)=0$ 
%but from Weyl-Heisenberg's inequality  and the $L^2$ conservation law \eqref{Conservation Law} it follows that 
%\[0<\|u_0\|^2=\|u(\cdot, t)\|_{L^2}^2\le \frac12 \bigl \||x|u(\cdot, t) \bigr\|_{L^2}
%\|\nabla u (\cdot,t)\|_{L^2}.\]
%Hence,
%\[\lim\limits_{t\nearrow t^*}\|\nabla u (\cdot,t)\|_{L^2}=+\infty.\]

\smallskip
Also, from \eqref{Energy-2}, H\"older's inequality and the Sobolev inequality in four dimensions,
we obtain the estimate 
\begin{equation*}
\begin{split}
\int_{\R^4} |\nabla u|^2 + v^2 dx &=E(t) - 2 \int_{\R^4} v |u|^2 dx \\
&\leq E_0 + 2 \int_{\R^4} |v| |u|^2 dx\\
&\leq E_0 + 2 \| \bar{u}\, v\|_{L^2_x} \|u\|_{L^2_x} \\
&\leq E_0 + 2 \|u\|_{L^4_x} \|v\|_{L^4_x} \|u_0\|_{L^2_x} \\
&\leq E_0 + 2\,c^2 \|\nabla u\|_{L^2_x} \|\nabla v\|_{L^2_x} \|u_0\|_{L^2_x} \\
&\leq E_0 + \frac12 \|\nabla u\|_{L^2_x}^2+  2\,c^4  \|\nabla v\|_{L^2_x}^2 \|u_0\|_{L^2_x}^2.
\end{split}
\end{equation*}
Therefore,
$$\frac12 \|\nabla u\|_{L^2_x}^2+\|v\|_{L^2_x}^2 \leq E_0 +  2\,c^4  \|\nabla v\|_{L^2_x}^2 \|u_0\|_{L^2_x}^2,$$
which, from \eqref{blowupfinal}, yields
\begin{equation}\label{chopp}
\lim\limits_{t\nearrow t^*}\|\nabla v (\cdot,t)\|_{L^2}=+\infty.
\end{equation}

\smallskip
Now, from the integral solution formula for $v$ \eqref{Sd-Integral-v}, we have for all $t <t^*$,
\begin{equation*}
\begin{split}
\|\nabla v(t)\|_{L^2} &\leq e^{-t/\mu} \|\nabla v_0 \|_{L^2}+ \tfrac{2}{\mu}\int_0^t\,
e^{-(t-\tau)/\mu} \|u(\tau)\|_{L^{\infty}} \|\nabla u(\tau)\|_{L^2}\,d\tau \\
& \leq \|\nabla v_0 \|_{L^2}+ 2(1-e^{-t/\mu}) \|u\|_{L^{\infty}_{[0,t]}L^{\infty}_x(\R^4)} \|\nabla 
u\|_{L^{\infty}_{[0,t]}L^2_x(\R^4)},
\end{split}
\end{equation*}
and therefore \eqref{chopp} implies that, either $\|u\|_{L^{\infty}_{[0,t^*]}L^{\infty}_x(\R^4)}=+\infty$
or $\|\nabla u\|_{L^{\infty}_{[0,t^*]}L^2_x(\R^4)}=+\infty$. The high regularity of the initial data
guarantees that the blowup of the mixed norms can only happen at $t=t^*$, i.e.
$\lim\limits_{t\nearrow t^*}\|u\|_{L^{\infty}_x(\R^4)}=+\infty$
or
$\lim\limits_{t\nearrow t^*}\|\nabla u\|_{L^2_x(\R^4)}=+\infty$. We will now see that, actually, both norms
explode at $t^*$.

\smallskip
In fact from
\eqref{Energy-1} and \eqref{Energy-2} the following inequality holds:
$$\int_{\R^n} |\nabla u|^2 + v^2 dx = E(t) - 2 \int_{\R^n} v |u|^2 dx \leq E_0 + 2 \int_{\R^n} |v| |u|^2 dx.$$
Using the estimate 
$$\int |v| |u|^2 dx \leq \|v\|_{L^1_x} \|u\|_{L^{\infty}_x}^2,$$
and the a priori bound \eqref{apriori_v} for the $L^1$ norm of $v$, we again obtain 
$\lim\limits_{t\nearrow t^*}\|u \|_{L^{\infty}_x(\R^4)}=+\infty$
when $\lim_{t\nearrow t^*} \|\nabla u\|_{L^2_x}=+\infty$, so that we conclude that 
$\|u \|_{L^{\infty}_x(\R^4)}$ always blows up at $t^*$.

\smallskip
On the other hand, from Sobolev's inequality, if $\lim\limits_{t\nearrow t^*}\|u \|_{L^{\infty}_x(\R^4)}=+\infty$ does happen,
then $\lim_{t\nearrow t^*} \| u\|_{H^s_x}=+\infty$ for $s>2$. A persistence of regularity argument,
analogous to Remark \ref{persistence}, but applied only to the $u$ component of the solution in the 
integro-differential formulation of the problem \eqref{SD-IDF}, implies that, for the higher regularity
norms of $u$ to blow up at $t^*$ it is necessary for that to happen to $s=1$ as well. Therefore we conclude
that, if $\lim\limits_{t\nearrow t^*}\|u \|_{L^{\infty}_x(\R^4)}=+\infty$ holds, then 
$\lim_{t\nearrow t^*} \| u\|_{H^1_x}=+\infty$ also holds which, due to the conservation of the $L^2$ norm, finally
implies $\lim_{t\nearrow t^*} \|\nabla u\|_{L^2_x}=+\infty$. This concludes the proof. $\square$

\smallskip
\section{\bf{Global Well-posedness for Small Regular Data in Four Dimensions}} 
\label{GWP-4D}

We end this work by using the persistence property, in the same way as was used
in the last paragraph of the proof of Theorem~\ref{blowup-theorem}, to show the existence of global solutions in
$H^1(\R^4) \times  H^1(\R^4)$, for
\eqref{SD}, in the focusing regime ($\lambda =-1$), and small regular intial data.  

\smallskip
\subsection{Proof of Proposition \ref{GWP-Small-4D}}
	
For data $(u_0, v_0)$ belonging  to the space $H^s(\R^4)\times H^s(\R^4)$  with $s>2$, consider the corres\-ponding $H^1(\R^4)\times H^1(\R^4)$ solution  $\bigl(u(\cdot, t),v(\cdot,t)\bigl)$ of (\ref{SD}) established in Theorem \ref{LWP-4D} and defined on its maximal positive time interval $\bigl[0, T_*\bigr).$ Using the pseudo-energy (\ref{Energy-2}) and combining the H\"older, Sobolev and Young inequalities we have
\begin{equation}
\begin{split}\label{proof-gwp-small-4D-1}
\|\nabla u(\cdot, t)\|_{L^2}^2 + \|v(\cdot, t)\|_{L^2}^2&\le E_0 -2\int_{\R^4}v(\cdot, t)|u(\cdot, t)|^2dx\\
&\le E_0 + 2 \|v(\cdot, t)\|_{L^2}\|u(\cdot, t)\|^2_{L^4}\\
&\le E_0 + 2 c_4^2 \|v(\cdot, t)\|_{L^2}\|\nabla u(\cdot, t)\|^2_{L^2}\\
&\le E_0 + \|v(\cdot, t)\|^2_{L^2} +  c_4^4\|\nabla u(\cdot, t)\|^4_{L^2},
\end{split}
\end{equation}
which yields the estimate
\begin{equation}\label{proof-gwp-small-4D-2}
\|\nabla u(\cdot, t)\|_{L^2}^2  \le E_0 +  c_4^4\|\nabla u(\cdot, t)\|^4_{L^2}. 
\end{equation}

\smallskip
Now, we define the non-negative continuous function
$x(t):=\|\nabla u(\cdot, t)\|_{L^2}^2,$ for all $0\le t < T_*$, that satisfies the inequality
\begin{equation}\label{proof-gwp-small-4D-3}
y(t):=x(t)-c_4^4x^2(t) \le E_0. 
\end{equation}

\smallskip
 From \eqref{proof-gwp-small-4D-3}, using classical bootstrap arguments, we can conclude that if the initial data verifies 
 the conditions:  
 \begin{equation}\label{proof-gwp-small-4D-4}
 0\leq E_0 < \beta:=\frac{1}{4c_4^4}\quad \text{and}\quad x(0)=\|\nabla u_0\|_{L^2}^2 \le \gamma_0:=\tfrac{1-\sqrt{1-4c_4^4E_0}}{2c_4^4},
 \end{equation}
 \smallskip 
 then (see Figure \ref{Fig-Bootstrap-4d}) we have
 \begin{equation}\label{proof-gwp-small-4D-5}
 \|\nabla u(\cdot, t)\|_{L^2}\le \sqrt{\gamma_0},\quad   \text{for all}\quad  t\in [0, T_*).
 \end{equation} 

\smallskip
\begin{figure}[ht]
	\begin{center}
		\begin{tikzpicture}[domain=0:5]
		\draw[->] (-0.5,0)--(5.5,0) node[below] {$x$};
		\draw[->] (0,-1)--(0,5) node[right] {$y$}; 
		\draw[thick, domain=0:4.3] plot (\x, 4*\x - \x^2) node at (4.2,3) {{\small $y=x-c_4^4x^2$}};
		\draw[thin, color=red, dashed] (-1,4.02)--(2.1, 4.02);
		\node at (2.,4){{\small $\bullet$}};
		\node at (-1.3,4.02){{\small $\frac{1}{4c_4^4}$}};
		\draw[thin, color=red, dashed] (0,3.8)--(5, 3.8);
		\node at (-0.3,3.7){{\small $E_0$}};
		\node at (2,0){{\small $\bullet$}};
		\node at (2.05,-0.35){\small{$\frac{1}{2c_4^4}$}};
		\node at (4,0){{\small $\bullet$}};
		\node at (3.8,-0.35){\small{$\frac{1}{c_4^4}$}};
		\draw[thin, dashed] (2,0)--(2, 4);
		\node at (1.6,0){{\small $\bullet$}};
		\draw[thin, dashed] (1.6,0)--(1.6, 3.8);
		\node at (1.55,-0.3){{\small $\gamma_0$}};
		\draw[line width =1.75](0,0)--(1.6, 0);
		\end{tikzpicture}
	\end{center}
	\vspace{-0.3cm}
	\caption{}
	\label{Fig-Bootstrap-4d}
\end{figure}

\smallskip
Besides, this a priori uniform control in time for $\|\nabla u(\cdot, t)\|_{L^2}$, under smallness conditions for the intital data, can be used to control 
$\|v(\cdot, t)\|_{L^2}$ similarly to the calculations performed in the proof of Theorem \ref{GWP-Small-3D}. However, we opted here to provide another approach for bounding
$\|v(\cdot, t)\|_{L^2}$ using \eqref{equation-v}, which combined with Sobolev embedding  and \eqref{proof-gwp-small-4D-5} gives us
\begin{equation*}\begin{split}
\|v(\cdot, t)\|_{L^2}&\le e^{-t/\mu}\|v_0\|_{L^2}+\frac{1}{\mu}\int_0^te^{-(t-\tau)/\mu}\|u(\cdot, \tau)\|_{L^4}^2d\tau\\
&\le e^{-t/\mu}\|v_0\|_{L^2} +c_4^2\frac{1}{\mu}\int_0^te^{-(t-\tau)/\mu}\|\nabla u(\cdot, \tau)\|_{L^2}^2d\tau\\
&\le e^{-t/\mu}\|v_0\|_{L^2} +c_4^2\gamma_0\frac{1}{\mu}\int_0^te^{-(t-\tau)/\mu}d\tau\\
&=e^{-t/\mu}\|v_0\|_{L^2} +c_4^2\gamma_0(1-e^{-t/\mu})\\
&\le \|v_0\|_{L^2} +c_4^2\gamma_0,
\end{split}\end{equation*}
for all $t\in [0, T_*)$.

\smallskip
It only remains to understand what happens with the growth of the $L^2$-norm of $\nabla v$. In this sense, as pointed out in last arguments in the proof of Theorem \ref*{blowup-theorem},
we know that if $\|\nabla v(\cdot, t)\|_{L^2}$ blows up in finite time, this would imply the existence of $t^*\leq T_*$ such that  $\lim\limits_{t\nearrow t^*}\|u \|_{L^{\infty}_x(\R^4)}=+\infty$. However, the
regularity of the initial data ($s>2$) and the Sobolev embedding, with the persistency property, would lead to a contradiction with \eqref{proof-gwp-small-4D-5}. Therefore, the boundedness of the $H^1\times H^1$ norm in the maximal time interval $[0, T^*)$ implies that there can be no blow-up of this norm in finite time, which gives us  $T^*=+\infty$. This completes the proof. $\square$

% % % % % % % % % % % % % % % % % % % % % % % % % % % % % % % % % % % % % % % % % % % % % % % % % % % % % %
\smallskip
\subsection*{Acknowledgments}
J. Drumond Silva would like to thank the kind hospitality of IMPA, Instituto de Matem\'atica Pura e Aplicada, and of the Institute of Mathematics at the Federal University of Rio de Janeiro, where part of this work was developed. Both
authors thank S. Correia and F. Oliveira for crucial observations and comments with respect to a previous version
of the manuscript.

% % % % % % % % % % % % % % % % % % % % % % % % % % % % % % % % % % % % % % % % % % % % % % % % % % % % % %


\begin{thebibliography}{99}

\bibitem{Besse-Bidegaray}{C. Besse and B. Bid\'egaray,}
{\it  Numerical study of self-Focusing solutions to the Schr\"odinger-Debye system,}
{ESAIM: M2AN, {\bf 35}  (2001), 35--55.}

\bibitem{Bidegaray1}{B. Bid\'egaray,}
{\it On the Cauchy problem for systems occurring in nonlinear optics,}
{Adv. Diff. Equat., {\bf 3} (1998), 473--496.}

\bibitem{Bidegaray2}{B. Bid\'egaray,}
{\it  The Cauchy problem for Schr\"odinger-Debye equations,}
{Math. Models Methods Appl. Sci., {\bf 10} (2000), 307--315.}

\bibitem{Carvajal}{X. Carvajal and P. Gamboa,}
{\it  Global well-posedness for the critical Schr\"odinger-Debye system,}
{Dynamics of PDE, {\bf 11} No.3 (2014), 251--268.}

\bibitem{Cazenave-Book}{T. Cazenave,}
{\it Semilinear Schr\"odinger equations,}
{Courant Lecture Notes in Mathematics, vol. 10, New York: New York University Courant Institute of Mathematical Sciences, AMS, 2003.}

\bibitem{Cazenave-Critical-1}{T. Cazenave and F. B. Weissler,}
{\it Some remarks on the nonlinear Schr\"odinger equation in the critical case,}
{Lecture Notes in Math., {\bf 1394} (1989), 18--29.}

\bibitem{Cazenave-Critical-2}{T. Cazenave and F. B. Weissler,}
{\it The Cauchy problem for the critical nonlinear Schr\"odinger equation in $H^s$,}
{Nonlinear Anal. TMA., {\bf 14} (1990), 807--836.}

\bibitem{Imethod}{J. Colliander, M. Keel, G. Staffilani, H. Takaoka and T. Tao,}{\it Global well-posedness for KdV in Sobolev spaces of negative index,} {Electron. J. Differential Equations, {\bf 26} (2001), 1--7.}

\bibitem{global3}{J. Colliander, M. Keel, G. Staffilani, H. Takaoka and T. Tao,}{\it Global existence and scattering for rough
	solutions of a nonlinear Schr\"odinger equation on $\R^3$,} {Comm. Pure Appl. Math., {\bf 57} (2004), 987--1014.}


\bibitem{global2}{J. Colliander and T. Roy,}{\it Bootstrapped Morawetz estimates and resonant decomposition for low regularity
global solutions of cubic NLS on $\R^2$,}{Commun. Pure Appl. Anal., {\bf 10} (2011), 397--414.}



\bibitem {Corcho-Linares}{A. J. Corcho and F. Linares,} \textit{Well-posedness
for the Schr\"odinger-Debye Equation,} {Contemporary Mathematics, \textbf{362}
(2004), 113--131.}

\bibitem {Corcho-Matheus}{A. J. Corcho and C. Matheus,} \textit{Sharp bilinear
estimates and well-posedness for the 1-D Schr\"{o}dinger-Debye system,}
{Differential and Integral Equations, \textbf{22} (2009), 357--391.}

\bibitem {Corcho-Oliveira-Silva}{A. J. Corcho, F. Oliveira and J. D. Silva} \textit{Local
and Global well-posedness for the critical Schr\"{o}dinger-Debye system,}
{Proceedings of the AMS, {\bf 141} (2013),  3485--3499.}

\bibitem{Fibich-Papanicolau}{G. Fibich and G. C. Papanicolau,}
{\it Self-focusing in the perturbed and unperturbed nonlinear
Schr\"odinger equation in critical dimension,}
{SIAM J. Appl. Math., {\bf 60} (1999), 183--240.}

\bibitem{Ginibre}{J. Ginibre, Y. Tsutsumi and G. Velo,}
{\it On the Cauchy problem for the Zakharov system,}
{J. Funct. Anal., {\bf 151} (1997), 384--436.}

\bibitem{Ginibrevelo}{J. Ginibre and G. Velo,}
{\it On a class of nonlinear Schr\"{o}dinger equations,}
{J. Funct. Anal., {\bf 32} (1979), 1-32, 33--72.}


\bibitem{Keel-Tao}{M. Keel and T. Tao,}
{\it Endpoint Strichartz estimates,}
{American Journal of Mathematics, {\bf 120} (1998), 955--980.}

\bibitem{KPV0-a}{C. E. Kenig, G. Ponce and L. Vega,}
{\it Small solutions to nonlinear Schr\"odinger equations,}
{Ann. Inst. Henri Poincar\'e, Analyse nonlin\'eaire, {\bf 10} (1993), 255--288.}


\bibitem{KPV0-b}{C. E. Kenig, G. Ponce and L. Vega,}
{\it Well-posedness and scattering results for the generalized Korteweg-de Vries equation
via the contraction principle,}
{Comm. Pure Appl. Math., {\bf 46} (1993), 527--620.}

\bibitem{Linares-Ponce}{F. Linares and G. Ponce,}
{\it ``Introduction to nonlinear dispersive equations'',}
{Springer, 2009.}


\bibitem{Merle}{F. Merle,}
{\it Blow-up results of viriel  type for Zakharov equations,}
{Comm. Math. Phys., {\bf 175} (1996), 433--455.}


\bibitem{Newell}{A. C. Newell and J. V. Moloney,}
{\it ``Nonlinear optics'',}
{Addison-Wesley, 1992.}

\bibitem{Strauss}{W. A. Strauss,}
{\it Existence of solitary waves in higher dimensions,}
{Comm. Math. Phys., {\bf 55} (1977), 149--162.}


\bibitem{Y-Tsutsumi}{Y. Tsutsumi,}
{\it $L^2$-solutions for nonlinear Schr\"odinger equations and nonlinear groups,}
{Funkcial. Ekvac., {\bf 30} (1987), 115--125.}

\end{thebibliography}
\end{document}